\documentclass{gtart}

\def\ifplaintex{\expandafter\ifx\csname documentclass\endcsname\relax}

\ifplaintex 
\else
\fi

\expandafter\ifx\csname beginpicture\endcsname\relax
\expandafter\ifx\csname documentclass\endcsname\relax
\input pictex \else
\input prepictex \input pictex \input postpictex \fi\fi

\def\gt{{\mathsurround=0pt\it $\cal G\mskip-2mu$eometry \&\ 
$\cal T\!\!$opology}}        

\def\gtp{{\mathsurround=0pt\it $\cal G\mskip-2mu$eometry \&\ 
$\cal T\!\!$opology $\cal P\!$ublications}}  

\def\lognumber#1{\def\thelognumber{#1}}
\def\volumenumber#1{\def\thevolumenumber{#1}}
\def\papernumber#1{\def\thepapernumber{#1}}
\def\volumeyear#1{\def\thevolumeyear{#1}}

\def\pagenumbers#1#2{\def\startpage{#1}\def\finishpage{#2}}
\def\published#1{\def\publishdate{#1}}
\def\proposed#1{\def\theproposer{#1}}
\def\seconded#1{\def\theseconders{#1}}
\def\received#1{\def\receiveddate{#1}}
\def\revised#1{\def\reviseddate{#1}}
\def\accepted#1{\def\accepteddate{#1}}

\def\asciiaddress#1{\def\theasciiaddress{#1}}
\def\asciiemail#1{\def\theasciiemail{#1}}
\long\def\asciiabstract#1{\long\def\theasciiabstract{#1}}

\let\\\par\let\thelognumber\relax
\let\thevolumenumber\relax\let\thepapernumber\relax
\let\thevolumeyear\relax\let\thesamplenumber\relax\let\startpage\relax
\let\finishpage\relax\let\publishdate\relax\let\receiveddate\relax
\let\reviseddate\relax\let\accepteddate\relax\let\theasciititle\relax
\let\theasciiauthors\relax\let\theasciiaddress\relax
\let\theasciiabstract\relax
\let\theasciiemail\relax\let\theshortauthors\relax\let\theshorttitle\relax

\long\def\maketitlep{   

\count0=\startpage

\gt\hfill      
\beginpicture
\setcoordinatesystem units <0.33truein, 0.33truein> point at 2.2 0.9
\setplotsymbol ({$\cal G$})
\plotsymbolspacing=9truept
\circulararc 315 degrees from 0 1 center at 0 0
\setplotsymbol ({$\cal T$})
\circulararc 315 degrees from 1 -1 center at 1 0
\endpicture
\break
{\small\ifx\thesamplenumber\relax 
Volume \else Sample
\fi\thevolumenumber\ (\thevolumeyear)
\startpage--\finishpage\nl
Published: \publishdate}
\vglue 0.5truein plus 0.4fil minus 0.1truein

{\parskip=0pt\leftskip 0pt plus 1fil\def\\{\par\smallskip}{\ifplaintex\large
\else\Large\fi\bf\thetitle}\par\medskip}   

\vglue 0pt plus 0.1fil 

{\parskip=0pt\leftskip 0pt plus 1fil\def\\{\par}{\sc\theauthors}
\par\medskip}

\vglue 0pt plus 0.1fil 

{\small\parskip=0pt\let\newline\\
{\leftskip 0pt plus 1fil\def\\{\par}{\sl\theaddress}\par}
\expandafter\ifx\theemail\relax    
\relax\else\vglue 5pt plus 0.02fil minus 2pt\def\\{\stdspace{\rm 
and}\stdspace} 
\cl{Email:\stdspace\tt\theemail}\fi
\ifx\theurl\relax                  
\relax\else\vglue 5pt plus 0.02fil minus 2pt\def\\{\stdspace{\rm 
and}\stdspace}
\cl{URL:\stdspace\tt\theurl}\fi\par}

\vglue 7pt plus 0.3fil minus 3pt

{\bf Abstract}
\vglue 5pt plus 0.1fil minus 2pt

\theabstract

\vglue 7pt plus 0.3fil minus 3pt

{\bf AMS Classification numbers}\quad Primary:\quad \theprimaryclass

Secondary:\quad \thesecondaryclass

\vglue 5pt plus 0.3fil minus 2pt

{\bf Keywords}\quad \thekeywords

\vglue 10pt plus 0.5fil minus 5pt

{\small  Proposed: \theproposer\hfill Received: \receiveddate\nl
Seconded: \theseconders\hfill 
\ifx\reviseddate\relax                         
Accepted: \accepteddate                        
\else
Revised: \reviseddate                          
\fi}
\eject
}       

\let\maketitlepage\maketitlep
\let\maketitle\maketitlepage

\font\phead=cmsl9 scaled 950
\font=cmsl9 scaled 1050
\font\pnum=cmbx10 scaled 913
\font\lnum=cmbx10 
\font\pfoot=cmsl9 scaled 950
\font=cmsl9 scaled 1050
\ifplaintex
\headline{\vbox to 0pt{\vskip -4.5mm\line{\small\phead\ifnum
\count0=\startpage ISSN 1364-0380 (on line)
1465-3060 (printed) \hfill {\pnum\folio}\else\ifodd\count0\def\\{ }%
\ifx\theshorttitle\relax\thetitle\else\theshorttitle\fi\hfill{\pnum\folio}
\else\def\\{ and }{\pnum\folio}\hfill\ifx\theshortauthors\relax\theauthors
\else\theshortauthors\fi\fi\fi}\vss}}
\footline{\vbox to 0pt{\vglue 0mm\line{\small\pfoot\ifnum\count0=\startpage
\copyright\ \gtp\hfill\else
\gt, Volume \thevolumenumber\ (\thevolumeyear)\hfill\fi}\vss
}}
\else
\makeatletter
\def\@oddhead{{\small\lhead\ifnum\count0=\startpage ISSN 1364-0380 (on line)
1465-3060 (printed) \hfill {\lnum\number\count0}\else\ifodd\count0
\def\\{ }\ifx\theshorttitle\relax \thetitle \else\theshorttitle\fi\hfill
{\lnum\number\count0}\else\def\\{ and }{\lnum\number\count0}
\hfill\ifx\theshortauthors\relax 
\theauthors\else\theshortauthors\fi\fi\fi}}\def\@evenhead{\@oddhead}
\def\@oddfoot{\small\lfoot\ifnum\count0=\startpage\copyright\ \gtp\hfill\else
\gt, Volume \thevolumenumber\ (\thevolumeyear)\hfill\fi}
\def\@evenfoot{\@oddfoot}
\makeatother
\fi

\newwrite\gtoutfile
\long\gdef\makeheadfile{  
{\def\\{, }\def\s{ }
\immediate\openout\gtoutfile head.xxx
\immediate\write\gtoutfile{Proxy-for: \ifx\theasciiauthors\relax
\theauthors\else\theasciiauthors\fi\s<\ifx\theasciiemail\relax\theemail\else\theasciiemail\fi>}
\immediate\write\gtoutfile{\noexpand\\}
\immediate\write\gtoutfile{Authors: \ifx\theasciiauthors\relax
\theauthors\else\theasciiauthors\fi}
{\def\\{ }\immediate\write\gtoutfile{Title: \ifx\theasciititle\relax
\thetitle\else\theasciititle\fi}}
\immediate\write\gtoutfile{Subj-class: GT or SG or MG etc}
\immediate\write\gtoutfile{MSC-class: \theprimaryclass\ifx\thesecondaryclass\relax\else, \thesecondaryclass\fi}
\immediate\write\gtoutfile{Journal-ref: Geom. Topol. \thevolumenumber
(\thevolumeyear) \startpage-\finishpage}
\immediate\write\gtoutfile{Comments: Published by Geometry and Topology at}
\immediate\write\gtoutfile{\s\s http://www.maths.warwick.ac.uk/gt/GTVol\thevolumenumber/paper\thepapernumber.abs.html}
\immediate\write\gtoutfile{\noexpand\\}
\immediate\write\gtoutfile{}
\ifx\theasciiabstract\relax
\immediate\write\gtoutfile{\theabstract}\else
\immediate\write\gtoutfile{\theasciiabstract}\fi
\immediate\write\gtoutfile{}
\immediate\write\gtoutfile{\noexpand\\}
\immediate\write\gtoutfile{}
\immediate\closeout\gtoutfile}}  

\def\maketitlepage{\maketitlep\makeheadfile}
\let\maketitle\maketitlepage

\lognumber{335}
\volumenumber{8}\papernumber{2}\volumeyear{2004}
\pagenumbers{35}{76}
\received{30 June 2003}
\revised{18 December 2003}
\published{21 January 2004}
\accepted{16 January 2004}
\proposed{Ronald Stern}
\seconded{Robion Kirby, Tomasz Mrowka}

\usepackage{amsmath,amssymb}
\usepackage{amscd,graphicx,psfrag,verbatim}

\nocolon

\newcommand{\R}{\mathcal R}
\newcommand{\p}{\partial}
\newcommand{\A}{\mathcal A}
\newcommand{\B}{\mathcal B}
\newcommand{\C}{\mathcal C}
\newcommand{\M}{\mathcal M}
\newcommand{\G}{\mathcal G}
\newcommand{\F}{\mathcal F}
\newcommand{\Z}{\mathbb Z}
\newcommand{\Q}{\mathbb Q}
\renewcommand{\P}{\mathcal P}
\renewcommand{\S}{\mathcal S}
\renewcommand{\H}{\mathcal H}
\newcommand{\ind}{\operatorname{ind}}
\newcommand{\hol}{\operatorname{hol}}
\newcommand{\ep}{\varepsilon}
\renewcommand{\phi}{\varphi}
\newcommand{\id}{\operatorname{id}}

\newcommand{\tr}{\operatorname{tr}}
\newcommand{\Id}{\operatorname{Id}}

\newcommand{\su}{\mathfrak{su}}

\renewcommand{\sf}{\operatorname{sf}}
\newcommand{\lk}{\operatorname{lk}}
\newcommand{\ad}{\operatorname{ad}}

\newcommand{\arf}{\operatorname{arf}}
\newcommand{\Fix}{\operatorname{Fix}}
\newcommand{\coker}{\operatorname{coker}}
\newcommand{\sign}{\operatorname{sign}}
\newcommand{\diag}{\operatorname{diag}}
\newcommand{\Hess}{\operatorname{Hess}}

\theoremstyle{plain} \newtheorem{theorem}{Theorem}[section]
\newtheorem{lemma}[theorem]{Lemma}
\newtheorem{corollary}[theorem]{Corollary}
\newtheorem{proposition}[theorem]{Proposition}

\theoremstyle{definition} \newtheorem{remark}[theorem]{Remark}

\begin{document}

\title{Rohlin's invariant and gauge theory II. \\ Mapping tori}
\authors{Daniel Ruberman\\Nikolai Saveliev}
\address{Department of Mathematics, MS 050, Brandeis
University\\Waltham, MA 02454, USA}
\email{ruberman@brandeis.edu}
\secondaddress{Department of Mathematics, University of Miami\\PO 
Box 249085, Coral Gables, FL 33124, USA}
\secondemail{saveliev@math.miami.edu}
\asciiaddress{Department of Mathematics, MS 050, Brandeis
University\\Waltham, MA 02454, USA\\and\\Department of Mathematics, University of Miami\\PO 
Box 249085, Coral Gables, FL 33124, USA}
\asciiemail{ruberman@brandeis.edu, saveliev@math.miami.edu}

\begin{abstract}
This is the second in a series of papers studying the relationship between
Rohlin's theorem and gauge theory. We discuss an invariant of a homology
$S^1 \times S^3$ defined by Furuta and Ohta as an analogue of Casson's
invariant for homology 3--spheres. Our main result is a calculation of the
Furuta--Ohta invariant for the mapping torus of a finite-order diffeomorphism
of a homology sphere. The answer is the equivariant Casson invariant
(Collin--Saveliev 2001) if the action has fixed points, and a version of the
Boyer--Nicas (1990) invariant if the action is free. We deduce, for
finite-order mapping tori, the conjecture of Furuta and Ohta that their
invariant reduces mod 2 to the Rohlin invariant of a manifold carrying a
generator of the third homology group. Under some transversality 
assumptions, we show
that the Furuta--Ohta invariant coincides with the Lefschetz number of the
action on Floer homology. Comparing our two answers yields an example of a
diffeomorphism acting trivially on the representation variety but
non-trivially on Floer homology.
\end{abstract}

\asciiabstract{This is the second in a series of papers studying the
relationship between Rohlin's theorem and gauge theory. We discuss an
invariant of a homology S^1 cross S^3 defined by Furuta and Ohta as an
analogue of Casson's invariant for homology 3-spheres. Our main result
is a calculation of the Furuta-Ohta invariant for the mapping torus of
a finite-order diffeomorphism of a homology sphere. The answer is the
equivariant Casson invariant (Collin-Saveliev 2001) if the action has
fixed points, and a version of the Boyer-Nicas (1990) invariant if the
action is free. We deduce, for finite-order mapping tori, the
conjecture of Furuta and Ohta that their invariant reduces mod 2 to
the Rohlin invariant of a manifold carrying a generator of the third
homology group. Under some transversality assumptions, we show that
the Furuta-Ohta invariant coincides with the Lefschetz number of the
action on Floer homology. Comparing our two answers yields an example
of a diffeomorphism acting trivially on the representation variety but
non-trivially on Floer homology.}

\primaryclass{57R57}
\secondaryclass{57R58} 

\keywords{Casson invariant, Rohlin invariant, Floer homology}
\maketitlepage

\section{Introduction}
Let $X$ be a $\Z[\Z]$--homology $S^1\times S^3$, that is, a smooth closed
oriented 4--manifold such that $H_* (X;\Z) = H_* (S^1 \times S^3;\Z)$ and
$H_*(\tilde X;\Z) = H_* (S^3;\Z)$, where $\tilde X$ is the universal
abelian cover of $X$. Denote by $\M^*(X)$ the moduli space of
irreducible ASD connections on a trivial $SU(2)$--bundle over $X$. The
virtual dimension of $\M^* (X)$ is $-3(1 - b_1 + b_2^+)(X) = 0$. In
fact, $\M^*(X)$ coincides with the moduli space of irreducible flat
connections on $X$. According to \cite{donaldson:orientation} and 
\cite{furuta-ohta}, the space $\M^*(X)$ is compact and canonically 
oriented once an orientation on $H^1(X;\mathbb R) = \mathbb R$ is fixed. 
After a perturbation if necessary, it is a finite collection of 
non-degenerate points. A signed count of these points is an invariant of 
$X$ known as the Donaldson polynomial $D_0(X)$ of degree zero. Furuta 
and Ohta \cite{furuta-ohta} define\footnote{Actually they use $D_0(X)/2$; 
we divide by $4$ to make the definition compatible with their conjecture, 
stated in the next paragraph.} $\lambda_{FO}(X) = D_0(X)/4$.  We will 
refer to $\lambda_{FO}(X)$ as the Furuta--Ohta invariant.

The Furuta--Ohta invariant may be regarded as a $4$--dimensional version 
of Casson's invariant.  Associated to $X$, there is also a Rohlin-type 
invariant, defined as the usual Rohlin invariant of a $3$--manifold 
carrying a generator of the third homology. (The definition of this 
invariant given in~\cite{furuta-ohta} is somewhat different from the one 
we use, which comes from~\cite{ruberman:ds}, but the two definitions can 
be readily shown to agree.)  Furuta and Ohta \cite{furuta-ohta} 
conjectured that the modulo 2 reduction of $\lambda_{FO} (X)$ equals the 
Rohlin invariant of $X$, by analogy with Casson's result that his 
invariant gives the Rohlin invariant of a homology $3$--sphere.

The goal of this paper is to calculate the Furuta--Ohta invariant in
the special case when
\[
X = [0,1]\times\Sigma\;/\;(0,x)\sim (1,\tau(x))
\]
is the mapping torus of a finite order orientation preserving
diffeomorphism $\tau\co 
\Sigma\to \Sigma$ of an integral homology sphere $\Sigma$.  The mapping
torus $X$ is a smooth 4--manifold oriented by the volume form $dt \wedge
\mathrm{vol}_\Sigma$. It is obviously a $\Z[\Z]$--homology $S^1\times
S^3$.  We show that the Furuta--Ohta invariant of $X$ is the equivariant
Casson invariant of the pair $(\Sigma,\tau)$.

More precisely, let $\R^* (\Sigma)$ be the space of irreducible flat
connections in the bundle $\Sigma\times SU(2)$ modulo gauge equivalence,
and let $\tau^*\co  \R^*(\Sigma)\to \R^* (\Sigma)$ be the map induced by
pulling back connections. The fixed
point set of $\tau^*$ will be denoted by $\R^{\tau}(\Sigma)$. After
perturbation if necessary, the space $\R^{\tau}(\Sigma)$ consists of
finitely many non-degenerate points which can be counted with signs
to obtain the equivariant Casson invariant $\lambda^{\tau}(\Sigma)$. A
rigorous definition of $\lambda^{\tau}(\Sigma)$ for $\tau$ having fixed
points can be found in \cite{collin-saveliev:casson}.  The definition of
$\lambda^{\tau}(\Sigma)$
for a fixed point free $\tau$ is given in this paper; it is related
to the Boyer--Nicas invariant~\cite{boyer-nicas,boyer-lines} of
$\Sigma/\tau$. Note that, if $\tau = \id$, the equivariant Casson invariant
coincides with the regular Casson invariant $\lambda(\Sigma)$.
\begin{theorem}\label{T:finite}
Let $\Sigma$ be an integral homology sphere and let $\tau\co \Sigma \to
\Sigma$ be
an orientation preserving diffeomorphism of finite order. If $X$ is the
mapping torus of $\tau$ then
\[
\lambda_{FO} (X) =  \lambda^{\tau} (\Sigma).
\]
\end{theorem}
We further show how to express $\lambda^{\tau} (\Sigma)$ in terms of the
regular Casson invariant and certain classical knot invariants, compare
with~\cite{collin-saveliev:casson}. Once such an explicit formula is in
place, we will prove in Section~\ref{S:rohlin} the following result, 
verifying the conjecture of Furuta and Ohta for the mapping tori of finite 
order diffeomorphisms.

\begin{theorem}\label{T:reduction}
Let $\tau\co  \Sigma \to \Sigma$ be an orientation preserving diffeomorphism
of finite order then the modulo 2 reduction of $\lambda^{\tau}(\Sigma)$
equals the Rohlin invariant of $\Sigma$.
\end{theorem}

An alternate approach to the invariant $\lambda_{FO}$ of a mapping
torus is via the `TQFT' view of Donaldson--Floer theory
\cite{donaldson:floer}.  The diffeomorphism $\tau$ induces an
automorphism on the instanton Floer homology of $\Sigma$.  Under the
assumption that the representation variety $\R^* (\Sigma)$ is non-degenerate,
we prove in Theorem~\ref{T:lefschetz} that the Furuta--Ohta invariant of $X$
equals half the Lefschetz number of this automorphism.  We conjecture that
this is the case in general, but the lack of equivariant perturbations is a
non-trivial obstacle to the proof.  A comparison of this result with
Theorem~\ref{T:finite} gives rise to the surprising phenomenon that a
map $\tau$ may act by the identity on $\R^* (\Sigma)$, but still act
non-trivially on the Floer homology.  Examples are given in
Section~\ref{S:examples}.

The basic idea behind Theorem~\ref{T:finite} is that for mapping tori,
$\lambda_{FO}$ counts fixed points of $\tau^*$, with multiplicity $2$
coming from the choice of holonomy in the circle direction.  Of
course, this must properly take
account of the signs with which flat connections are counted, and of
the perturbations used in the two theories.

The authors are grateful to Fred Diamond, Chris Herald, Chuck
Livingston, Tom Mrowka, and Liviu Nicolaescu for useful remarks and
for sharing their expertise.  The first author was partially supported
by NSF Grants 9971802 and 0204386. The second author was partially
supported by NSF Grant 0305946.


\section{Some equivariant gauge theory}\label{S:review}
Let $\A(\Sigma)$ be the affine space of connections in a trivialized
$SU(2)$--bundle $P$ over $\Sigma$, and let $\A^*(\Sigma)$ be the subset
consisting of irreducible connections. Any endomorphism $\tilde\tau\co 
P \to P$ which lifts $\tau\co  \Sigma \to \Sigma$ induces an action on
connections by pull back. For any two lifts, $\tilde\tau_1$ and
$\tilde\tau_2$, there obviously exists a gauge transformation $g\co  P \to
P$ such that $\tilde\tau_2 = \tilde\tau_1\cdot g$. This observation
shows that we have a well defined action $\tau^*$ on $\B^*(\Sigma) =
\A^*(\Sigma)/\G(\Sigma)$. The fixed point set of $\tau^*$ will be
called $\B^{\tau}(\Sigma)$.


\subsection{Decomposing $\B^{\tau}(\Sigma)$}
Let $\alpha$ be a connection on $P$ whose gauge equivalence class belongs
to $\B^{\tau} (\Sigma)$. Then there is a lift $\tilde\tau\co  P \to P$ such
that $\tilde\tau^*\alpha = \alpha$. Since $\alpha$ is irreducible, the
lift $\tilde\tau$ is defined uniquely up to a sign. Moreover,
$\tilde\tau^n$ is an endomorphism of $P$ lifting the identity map hence
$(\tilde\tau^n)^* \alpha = \alpha$ implies that $\tilde\tau^n = \pm 1$.
This allows one to decompose $\B^{\tau}(\Sigma)$ into a disjoint
union~\cite{braam-matic:smith,wang:involution}
\begin{equation}\label{E:decomp1}
\B^{\tau}(\Sigma) = \bigsqcup_{[\tilde\tau]}\; \B^{\tilde\tau}(\Sigma)
\end{equation}
where the $\tilde\tau$ are lifts of $\tau$ such that $\tilde\tau^n = \pm 1$.
The equivalence relation among the lifts $\tilde\tau$ is that $\tilde\tau_1
\sim \tilde\tau_2$ if and only if $\tilde\tau_2 = \pm\, g\cdot \tilde\tau_1
\cdot g^{-1}$ for some gauge transformation $g\co  P\to P$.

The spaces $\B^{\tilde\tau}(\Sigma)$ can be described as follows. For a
fixed lift $\tilde\tau$, let $\A^{\tilde\tau}(\Sigma)$ be the subset of
$\A^*(\Sigma)$ consisting of irreducible connections $A$ such that
$\tilde\tau^* A = A$. Define $\G^{\tilde\tau}(\Sigma) = \{\,g\in\G(\Sigma)
\; |\; g\tilde\tau = \pm \tilde\tau g\,\}$ then $\B^{\tilde\tau}(\Sigma) =
\A^{\tilde\tau}(\Sigma)/\G^{\tilde\tau}(\Sigma)$.

There is also a well defined action $\tau^*$ on the flat moduli space
$\R^*(\Sigma)$ whose fixed point set is $\R^{\tau}(\Sigma) = \B^{\tau}
(\Sigma)\,\cap\,\R^*(\Sigma)$. It splits as
\[
\R^{\tau}(\Sigma) = \bigsqcup_{[\tilde\tau]}\; \R^{\tilde\tau}(\Sigma).
\]
We wish to ramify the above splittings for later use. Note that any lift
$\tilde\tau\co  P \to P$ can be written in the base-fiber coordinates as
$\tilde\tau (x,y) = (\tau(x), \rho(x)\,y)$ for some $\rho\co  \Sigma \to
SU(2)$. The lift $\tilde\tau$ is said to be \emph{constant} if there
exists $u\in SU(2)$ such that $\rho(x) = u$ for all $x \in \Sigma$. The
rest of this section is devoted to proving the following result.

\begin{proposition}\label{P:lifts}
There are finitely many equivalence classes $[\tilde\tau]$ in the
decomposition (\ref{E:decomp1}), and each of them contains a constant lift.
\end{proposition}


\subsection{The case of non-empty $\Fix(\tau)$}
Let us first suppose that $\Fix(\tau)$ is non-empty. Then $\Sigma/\tau$
is an integral homology sphere and the projection map $\Sigma \to
\Sigma/\tau$ is a branched covering with branch set a knot. Let
$\tilde\tau (x,y) = (\tau(x), \rho(x)\,y)$ for some $\rho\co  \Sigma \to
SU(2)$. If $x \in \Fix(\tau)$ then $\tilde\tau(x,y) = (x,\rho(x)\,y)$
and $\tilde\tau^n (x,y) = (x,\rho(x)^n\,y) = (x,\pm y)$. This implies
that $\rho(x)^n =\pm 1$ and, in particular, that $\tr\rho(x)$ can only
take finitely many distinct values. Since $\Fix(\tau)$ is connected,
$\rho(\Fix(\tau))$ has to belong to a single conjugacy class in $SU(2)$.
Replacing $\tilde\tau$ by an equivalent lift if necessary, we may assume
that $\rho(x) = u$ for all $x \in \Fix(\tau)$.

Let $u\co  P \to P$ be the constant lift $u(x,y) = (\tau(x), u\cdot y)$
and consider the $SO(3)$ orbifold bundles $P/\tilde\tau$ and $P/u$
over $\Sigma/\tau$. All such bundles are classified by the holonomy
around the singular set in $\Sigma/\tau$. Since this holonomy equals
$\ad(u)$ in both cases, the bundles $P/\tilde\tau$ and $P/u$ are
isomorphic. Take any isomorphism and pull it back to a gauge
transformation $g\co  P\to P$. It is clear that $\tilde\tau = \pm\, g
\cdot u\cdot g^{-1}$ and hence $\B^{\tilde\tau}(\Sigma) = \B^u
(\Sigma)$. Thus (\ref{E:decomp1}) becomes a finite decomposition
\begin{equation}\label{E:decomp2}
\B^{\tau} (\Sigma) = \bigsqcup_{|\tr u|}\; \B^u (\Sigma).
\end{equation}


\subsection{The case of empty $\Fix(\tau)$}\label{S:decomp3}
Now suppose that $\Fix(\tau)$ is empty. Then $\Sigma/\tilde\tau$ is a
homology lens space and the projection map $\Sigma\to\Sigma/\tilde\tau$
is a regular (unbranched) covering. The bundle $P$ gives rise to the
$SO(3)$--bundle $P/\tilde\tau$ on $\Sigma/\tilde\tau$. Such bundles are
classified by their Stiefel--Whitney class $w_2(P/\tilde\tau)$ so there
are two different bundles if $n$ is even, and just one if $n$ is odd.

Note that the trivial bundle can be realized as the quotient bundle of
the constant lift $u(x,y) = (\tau(x), y)$, and the non-trivial bundle
(in the case of even $n$) as the quotient bundle of the constant lift
$u(x,y) = (\tau(x), u\cdot y)$ where $u$ is any matrix in $SU(2)$ such
that $u^n = -1$. Any isomorphism between $P/\tilde\tau$ and $P/u$, for
either of the above lifts $u$, pulls back to a gauge transformation $g\co 
P \to P$ such that $\tilde\tau =\pm\,g\cdot u\cdot g^{-1}$. Therefore,
the equivalence classes of lifts $\tilde\tau$ are classified by $w_2
(P/\tilde\tau)$.


\section{The unperturbed case}
In this section, we prove Theorem \ref{T:finite} under the assumption that
$\R^{\tau} (\Sigma)$ is non-degenerate, leaving the degenerate case to
later sections.  We first establish a two-to-one correspondence between
the flat moduli spaces $\M^*(X)$ and $\R^{\tau}(\Sigma)$, then compare
the non-degeneracy conditions in the two settings. Finally, the
orientations of $\M^*(X)$ and $\R^{\tau}(\Sigma)$ are compared using
the concept of orientation transport, see \cite{nicolaescu:swbook}.


\subsection{Identifying flat moduli spaces}\label{S:id-gauge}
Let $X$ be the mapping torus of $\tau\co  \Sigma\to\Sigma$. Denote by $\B^*(X)
= \A^* (X)/\G (X)$ the moduli space of irreducible connections in a trivial
$SU(2)$ bundle on $X$, and by $\M^* (X) \subset \B^* (X)$ the respective
anti-self-dual moduli space.  Since the bundle is trivial, a standard
Chern--Weil argument implies that $\M^* (X)$ consists of flat
connections.   Let $i\co  \Sigma \to X$ be the embedding $i(x) = [0,x]$.

\begin{proposition}\label{P1}
The map $i^*\co  \M^* (X) \to \R^{\tau} (\Sigma)$ induced by pulling back
connections is a well defined two-to-one correspondence.
\end{proposition}

\begin{proof}
First note that, for any irreducible flat connection $A$ on $X$, its pull
back $i^* A$ is also irreducible: if $i^* A$ were reducible it would have
to be trivial which would obviously contradict the irreducibility of $A$.

Let $P$ be a trivialized $SU(2)$--bundle over $\Sigma$. Given a flat
connection $A$ over $X$, cut $X$ open along $i(\Sigma)$ and put $A$ into
temporal gauge over $[0,1] \times \Sigma$. We obtain a path $A(t)$ of
connections in $P$. Note that $A(0)$ and $A(1)$ need
not be equal but they certainly are gauge equivalent (via the holonomy
along the intervals $[0,1]\times\{x\}$\,). Thus $A(1) = \tilde\tau^*
A(0)$ for some bundle automorphism $\tilde\tau\co  P \to P$ lifting $\tau$.
In temporal gauge, the flatness equation $F(A) = 0$ takes the form
$dt\wedge A'(t) + F(A(t)) = 0$ hence $A(t) = \alpha$ is a constant path,
where $\alpha$ is a flat connection over $\Sigma$. Since $\tilde\tau^*
A(0) = A(1)$ we conclude that $\tilde\tau^*\alpha = \alpha$.

Note that, conversely, $A$ can be obtained from $\alpha$ by pulling
$\alpha$ back to $[0,1] \times \Sigma$ and identifying the ends via
$\tilde\tau$.

Now, we need to see how the above correspondence behaves with respect
to gauge transformations. Let us fix a lift $\tilde\tau$.
Suppose that $A$ and $A'$ are connections that are in temporal gauge
when pulled back to $[0,1] \times \Sigma$, and that they are
equivalent via a gauge transformation $g$.  It is straightforward to
show that the restriction of $g\in \G(X)$ to $[0,1] \times \Sigma$
must be constant in $t$.
Since $g$ defines a gauge transformation over $X$, it also satisfies
the boundary
condition $\tilde\tau g = g \tilde\tau$. This identifies $\M^* (X)$
with the space of irreducible flat connections in $P$ modulo the
index two subgroup of $\G^{\tilde\tau}(\Sigma)$ which consists of
gauge transformations $g\co  P \to P$ such that $\tilde\tau g = g
\tilde\tau$. Because of the irreducibility, this leads to a two-to-one
correspondence between $\M^*(X)$ and $\R^{\tilde\tau}(\Sigma)$.
\end{proof}

\begin{remark}\label{R:pi}
Note that the construction of a connection over $X$ via pulling an
equivariant connection $\alpha$ back to $[0,1] \times \Sigma$ and
then identifying the ends makes sense for all connections $\alpha\in
\A^{\tilde\tau}(\Sigma)$ and not just the flat ones.  We will denote
the respective map by $\pi\co  \A^{\tilde\tau} (\Sigma) \to \A(X)$.
\end{remark}


\subsection{The non-degeneracy condition}
The moduli space $\R^{\tau}(\Sigma)$ is called {\it non--degenerate} if the
equivariant cohomology groups $H^1_{\tau}(\Sigma;\ad\alpha)$ vanish for all
$\alpha\in\R^{\tau}(\Sigma)$, compare with \cite{collin-saveliev:casson}.
The moduli space $\M^*(X)$ is called {\it non-degenerate} if
$\coker(d_A^* \oplus d_A^+) = 0$ for all $A \in \M^*(X)$.

\begin{proposition}\label{P2}
The moduli space $\M^*(X)$ is non-degenerate if and only if $\R^{\tau}
(\Sigma)$ is non-degenerate.
\end{proposition}

\begin{proof}
Since the formal dimension of $\M^*(X)$ is zero, $\ind(d_A^*\oplus d_A^+)
= 0$, and proving that $\coker (d_A^*\oplus d_A^+) = 0$ is equivalent to
proving that $\ker(d_A^*\oplus d_A^+) = 0$. The connection $A$ is flat
and therefore the latter is equivalent to showing that $H^1 (X;\ad A) = 0$.
The group $H^1(X;\ad A)$ can be computed with the help of the Leray--Serre
spectral sequence applied to the fibration $X\to S^1$ with fiber $\Sigma$.
The $E_2$--term of this spectral sequence is
\[
E_2^{pq} = H^p(S^1,\H^q(\Sigma;\ad\alpha)),
\]
where $\alpha = i^* A$ and $\H^q(\Sigma;\ad\alpha)$ is the local coefficient
system associated with the fibration. The groups $E_2^{pq}$ vanish for all
$p\ge 2$ hence the spectral sequence collapses at the second term, and
\[
H^1(X;\ad A) = H^1(S^1,\H^0(\Sigma;\ad\alpha))\,\oplus\,H^0(S^1,\H^1(\Sigma;
\ad\alpha)).
\]
Since $\alpha$ is irreducible, $H^0(\Sigma;\ad\alpha) = 0$ and the first
summand in the above formula vanishes. The generator of $\pi_1(S^1)$ acts on
$H^1 (\Sigma;\ad\alpha)$ as $\tau^*\co  H^1 (\Sigma;\ad\alpha)\allowbreak \to
H^1 (\Sigma;\ad\alpha)$, therefore, $H^0 (S^1, \H^1 (\Sigma;\ad\alpha))$ is
the fixed point set of $\tau^*$, which is the equivariant cohomology
$H^1_{\tau}(\Sigma;\ad\alpha)$. Thus we conclude that $H^1(X;\ad A) =
H^1_{\tau}(\Sigma;\ad\alpha)$, which completes the proof.
\end{proof}


\subsection{Orientation transport}
Let $M$ be a smooth closed oriented Riemannian manifold and let $D\co  C^{\infty}
(\xi)\to C^{\infty}(\eta)$ be a first order elliptic operator such that
$\ind D = 0$. Given a smooth family $h$ of bundle isomorphisms $h_s\co  \xi
\to\eta$, $0\le s\le 1$, referred to as \emph{homotopy}, form a family
of elliptic operators $D_s\co  C^{\infty}(\xi)\to C^{\infty}(\eta)$ by the
rule $D_s = D + h_s$. All these operators have the same symbol; in
particular, $\ind D_s = 0$ for all $s\in [0,1]$.

Let us fix orientations on the lines $\det D_i = \det(\ker D_i\,\otimes\,
(\coker D_i)^*)$, $i = 0, 1$. The homotopy $h$ provides an isomorphism
$\psi\co  \det D_0 \to \det D_1$. We say that the {\it orientation transport}
along homotopy $h$ is $1$ if $\psi$ is orientation preserving, and $-1$
otherwise. We use notation $\ep(D_0,h,D_1) = \pm 1$ for the orientation
transport, to indicate its dependence on the choice of orientations and
the homotopy $h$.

Once the orientations on $\det D_i$, $i = 0, 1$, are fixed, the
orientation transport $\ep(D_0,h,D_1)$ only depends on the homotopy class
of $h$ rel $\{\,0,1\,\}$. It is given by the following formula, see
\cite[page 95]{nicolaescu:swbook}. Define the {\it resonance set} of the
homotopy $h$ as
\[
\mathcal Z_h = \{\,s\in [0,1]\;|\;\ker D_s\ne 0\,\}.
\]
For each $s\in [0,1]$ denote by $P_s$ the orthogonal projection onto
$\coker D_s$. Let $h'_s$ be the derivative of the bundle isomorphism $h_s\co 
\xi\to\eta$. The homotopy is called {\it regular} if its resonance set is
finite and, for any $s\in [0,1]$, the {\it resonance operator}
\[
\begin{CD}
R_s = P_s\circ h'_s\co  \ker D_s @>>> L^2(\eta) @>>> \coker D_s
\end{CD}
\]
is a linear isomorphism. Suppose $h$ is a regular homotopy and set $d_s =
\dim\ker D_s = \dim\coker D_s$. Then
\begin{equation}\label{eq}
\ep (D_0,h,D_1) = \sign R_0\cdot\sign R_1\cdot\prod_{s\in [0,1)} (-1)^{d_s},
\end{equation}
where $\sign R_i = \pm 1$, $i = 0, 1$, according to whether $\det R_i\in
\det(\ker D_i \otimes (\coker D_i)^*)$ is positive or negative.


\subsection{Orientation of $\R^*(\Sigma)$}\label{S:R}
We assume that $\R^{\tau}(\Sigma)$ is non-degenerate. For any point
$\alpha\in\R^{\tau}(\Sigma)$, its orientation is given by
\[
(-1)^{\sf^{\tau}(\theta,\alpha)}
\]
where $\sf^{\tau}(\theta,\alpha)$ is the modulo 2 (equivariant) spectral
flow. In fact, equivariant spectral flow is well defined modulo 4 when
$\tau\ne\id$ and modulo 8 when $\tau =\id$ but the modulo 2 spectral flow
will suffice for our purposes. The definition is as follows, compare with
\cite{collin-saveliev:casson}.

Fix a Riemannian metric on $\Sigma$ so that $\tau\co \Sigma\to\Sigma$ is an
isometry, and consider a trivialized $SU(2)$--bundle $P$ over $\Sigma$.
According to Proposition \ref{P:lifts}, there exists a constant lift $u\co 
P \to P$ of $\tau$ such that $u^*\alpha = \alpha$ and $u^n = \pm 1$.

Observe that $u^*\theta = \theta$ where $\theta$ is the product connection
on $P$, and choose a smooth path $\alpha(s)$, $0\le s\le 1$, of equivariant
connections on $P$ such that $\alpha(0) = \theta$ and $\alpha(1) = \alpha$.
The equivariance here means that $u^*\alpha(s) = \alpha(s)$ for all $s$; it
can be achieved by averaging because connections form an affine space.
Associated with $\alpha(s)$ is a path of self-adjoint Fredholm operators
$K^u_{\alpha(s)}$ obtained by restricting the operators
\begin{equation}\label{E:ka}
K_{\alpha(s)} = \left(
\begin{array}{cr}
0             &   d^*_{\alpha(s)} \\
d_{\alpha(s)} &   -*d_{\alpha(s)}
\end{array}
\right)
\end{equation}
onto the space of $u$--equivariant differential forms $(\Omega^0\,\oplus
\,\Omega^1)_u(\Sigma,\ad P)$, where $\ad P = P\times_{\ad}\su(2)$ is the
adjoint bundle of $P$ and $d_{\alpha(s)}$ is the covariant derivative.
By $u$--equivariant differential forms $\Omega^*_u$ we mean the
$(+1)$--eigenspace of the pull back operator $u^*\co  \Omega^* \to
\Omega^*$ induced by the lift $\ad(u)\co  \ad P\to \ad P$.

The one--parameter family of spectra of operators $K^u_{\alpha(s)}$ can
be viewed as a collection of spectral curves in the $(s,\lambda)$--plane
connecting the spectrum of $K^u_{\theta}$ with that of $K^u_{\alpha}$.
These curves are smooth, at least near zero. The (modulo 2) spectral flow
$\sf^{\tau}(\theta,\alpha)$ is the number of eigenvalues, counted with
multiplicities, which cross the $s$--axis plus the number of spectral
curves which start at zero and go down. An equivalent way to define
spectral flow is to consider the straight line connecting the points
$(0,-\delta)$ and $(1,\delta)$ where $\delta > 0$ is chosen smaller than
the absolute value of any non-zero eigenvalue of $K^u_{\theta}$ and
$K^u_{\alpha}$. Then the spectral flow is the number of eigenvalues,
counted with multiplicities, which cross this line. The spectral flow is
well defined; it only depends on $\alpha$ and not on the choice of
$\alpha(s)$.

\begin{proposition}\label{P:reg}
Let $\alpha \in \R^\tau(\Sigma)$ and let $\alpha(s)$ be a path
of equivariant connections such that $\alpha(0) = \theta$ and $\alpha(1)
= \alpha$. Let $h_s = K^u_{\alpha(s)} - K^u_{\theta}$ then, for any small
generic $\delta > 0$,
\begin{equation}\label{eq2}
(-1)^{\sf^{\tau}(\theta,\alpha)} = \ep(K^u_{\theta} + \delta, h,
K^u_{\alpha} + \delta).
\end{equation}
\end{proposition}

\begin{remark}
If $\delta > 0$ is sufficiently small then both kernel and cokernel of
$K^u_{\alpha(s)} + \delta$ vanish at $s = 0$ and $s = 1$. This provides
canonical orientations for $\det(K^u_{\theta} + \delta)$ and
$\det(K^u_{\alpha} + \delta)$, which are implicit in (\ref{eq2}).
\end{remark}

\begin{proof}
The resonance set $\mathcal Z_h$ of $h_s$  consists of the values of $s$
at which the spectral curves of $K^u_{\alpha(s)} + \delta$ intersect the
$s$--axis. These are exactly the points where the spectral curves of
$K^u_{\alpha(s)}$ intersect the horizontal line $\lambda = - \delta$.
According to Sard's theorem, these intersections are transversal for a
generic $\delta > 0$. Moreover, the modulo 2 count of these intersection
points equals $\sf^{\tau}(\theta,\alpha)$ (due to the fact that $\alpha$
is non-degenerate and hence $\ker K^u_{\alpha} = 0$).

Now suppose that $s_0\in\mathcal Z_h$ and consider a smooth family
$\phi_s$ such that
\[
(K^u_{\alpha(s)} + \delta)\,\phi_s = \lambda_s\,\phi_s\quad\text{and}
\quad \lambda_{s_0}=0.
\]
Differentiating with respect to $s$ and setting $s = s_0$, we obtain
\[
h'_{s_0}\,\phi_{s_0} =
\lambda'_{s_0}\,\phi_{s_0} - (K^u_{\alpha(s_0)}+\delta)\,\phi'_{s_0}.
\]
The differential form $(K^u_{\alpha(s_0)} + \delta)\,\phi'_{s_0}$ is
orthogonal to $\coker(K^u_{\alpha(s_0)}+\delta)$ hence $\phi_{s_0}$ is
an eigenform of the resonance operator
\[
R_{s_0}\co  \ker (K^u_{\alpha(s_0)}+\delta) \to \coker (K^u_{\alpha(s_0)}+
\delta)
\]
with eigenvalue $\lambda'_{s_0}$. The path $\alpha(s)$ was chosen so that
$\lambda'_{s_0}\ne 0$ hence we can conclude that the resonance operator
$R_{s_0}$ is a linear isomorphism. Since this is true for all $s_0\in
\mathcal Z_h$, the homotopy $h_s$ is regular. Since $\ker (K^u_{\theta} +
\delta) = \ker (K^u_{\alpha} + \delta) = 0$ we conclude that $\sign R_0 =
\sign R_1 = 1$. The result now follows from (\ref{eq}).
\end{proof}


\subsection{Orientation of $\M^*(X)$}\label{S:orient4}
The mapping torus $X$ gets a canonical Riemannian metric once $\Sigma$ is
endowed with a Riemannian metric such that $\tau\co  \Sigma\to \Sigma$ is an
isometry. Fix a constant lift $u\co  P \to P$ and denote by $\bar P$ the
(trivial) $SU(2)$--bundle over $X$ obtained by first pulling $P$ back to
$[0,1]\times\Sigma$ and then identifying the ends via $u$. Given a
connection $A$ in $\bar P$, consider the respective ASD operator
\begin{equation}\label{E:da}
D_A = d_A^*\,\oplus\,d_A^+\co  \Omega^1(X,\ad\bar P) \to (\Omega^0\,\oplus\,
\Omega^2_+)(X,\ad \bar P).
\end{equation}
Observe for later use that the mapping torus structure on $X$ can be
used to identify both spaces $\Omega^1 (X,\ad\bar P)$ and $(\Omega^0\,
\oplus\,\Omega^2_+) (X,\ad\bar P)$ with the space $(\Omega^0\,\oplus\,
\Omega^1)^u (\mathbb R \times \Sigma,\ad P)$ consisting of $\su(2)$--valued
0-- and 1--forms $\omega(t,x)$ on $\Sigma$ depending on the parameter $t$
in a periodic fashion, $\omega (t+1,x) = u^*\omega (t,x)$. Here, $u$ can
be any lift of $\tau\co  \Sigma \to \Sigma$ to the bundle $P$. Note that any
periodic map $\mathbb R \to (\Omega^0\,\oplus\,\Omega^1)_u (\Sigma,\ad P)$
gives rise to an element of $(\Omega^0\,\oplus\,\Omega^1)^u(\mathbb R\times
\Sigma,\ad P)$. A form constructed in this way has the property that
$\omega(t,x) = u^* \omega
(t,x)$; not all elements of $(\Omega^0\,\oplus\,\Omega^1)^u(\mathbb R\times
\Sigma;\ad P)$ are of this kind. For any lift $u$ such that $u^*A = A$, the
operator $D_A$ can be viewed as
\[
D_A\co  (\Omega^0\,\oplus\,\Omega^1)^u (\mathbb R \times \Sigma,\ad P) \to
(\Omega^0\,\oplus\,\Omega^1)^u (\mathbb R \times \Sigma,\ad P).
\]
Let $\lambda_X$ be the determinant bundle of $D_A$ over $\B(X)$. This is
a real line bundle with the property that, over $\M^*(X) \subset \B(X)$,
the restriction of $\lambda_X$ is isomorphic to the orientation bundle of
$\M^*(X)$. According to \cite{donaldson:orientation}, the bundle $\lambda_X$
is trivial over $\B(X)$. Therefore, a choice of trivialization of $\lambda_X$
(given by a \emph{homology orientation}, that is, an orientation of
$H^1(X;\mathbb R) \oplus H^2_+(X;\mathbb R) = H^1 (X;\mathbb R)$\,) fixes an
orientation on $\M^*(X)$.

Since $\M^* (X)$ is non-degenerate, see Proposition \ref{P2}, it consists
of finitely many points. The above orientation convention translates into
orienting $[A]\in\M^*(X)$ by $\ep(D_{\theta},H,D_A)$, once orientations of
$\det D_A$ and $\det D_{\theta}$ are fixed. The former has a canonical
orientation because $\ker D_A = \coker D_A = 0$ due to the non-degeneracy.
The isomorphisms $\ker D_{\theta}= H^1 (X;\ad\theta)$ and $\coker
D_{\theta} = H^0(X;\ad\theta)$ provide a canonical orientation for
$\det D_{\theta}$ once we fix the orientation of the base of the fibration
$X\to S^1$.


\subsection{Proof of Theorem \ref{T:finite}}
We still assume that $\R^{\tau}(\Sigma)$ is non-degenera\-te. In this case,
$\lambda^{\tau}(\Sigma)$ is defined as half the signed count of points in
$\R^{\tau}(\Sigma)$ (see Section \ref{D:eq}), and $\lambda_{FO}(X)$ as one
fourth the signed count of points in $\M^*(X)$ (see Section \ref{D:fo}).
According to Proposition \ref{P1}, the restriction map $i^*\co  \M^*(X) \to
\R^{\tau} (\Sigma)$ is two-to-one, therefore, to prove the theorem, it is
sufficient to show that $i^*$ is orientation preserving.

Let $\alpha$ be an equivariant flat connection on $\Sigma$ and choose a
constant lift $u\co  P \to P$ of the diffeomorphism $\tau$ such that $\alpha
= u^*\alpha$ and $u^n = \pm 1$. By averaging, choose a path $\alpha(s)$,
$0\le s\le 1$, connecting $\theta$ to $\alpha$ such that $\alpha(s) =
u^*\alpha(s)$ for all $s$. It gives rise to a path of connections $A(s)$
on the product $[0,1]\times\Sigma$ given by the formula $A(s)(t,x) =
\alpha(s)(x)$, $x\in \Sigma$, compare with Section \ref{S:id-gauge}. These
define a path of connections on the mapping torus, called again $A(s)$,
because $A(s)(1,x) = u^* A(s)(0,x)$.

The following proposition is the main step in comparing the orientations.

\begin{proposition}\label{P:delta}
Let $A\in \M^*(X)$ and $\alpha = i^*A\in \R^{\tau}(\Sigma)$ then, for the
homotopy $H_s = D_{A(s)} - D_{\theta}$ and a small generic $\delta > 0$,
\[
\ep(D_{\theta} - \delta, H, D_A - \delta) = \ep(K^u_{\theta} + \delta, h,
K^u_{\alpha} + \delta).
\]
\end{proposition}

\begin{proof}
We begin by identifying $\ker(D_{A(s)} - \delta)$ with $\ker(K^u_{\alpha(s)}
+ \delta)$. Let us fix an $s$, and let $\{\,\psi_k\,\}$ be a basis of
eigenforms for the operator $K_{\alpha(s)}$,
\[
K_{\alpha(s)}\,\psi_k = \lambda_k\,\psi_k.
\]
Any differential form $\omega$ on $X$ can then be written in this
basis as
\[
\omega(t,x) = \sum_k\,a_k(t)\,\psi_k(x).
\]
The standard calculation shows that $D_{A(s)} = \p/\p t - K_{\alpha(s)}$
hence the equation $(D_{A(s)} - \delta)\omega = 0$ is equivalent to
\[
\sum_k\,(a'_k - (\lambda_k + \delta)\,a_k)\,\psi_k = 0,
\]
or $a'_k = (\lambda_k + \delta)\,a_k$ for every $k$. Therefore, $a_k(t) =
a_k(0)\,e^{(\lambda_k + \delta)\,t}$ and
\[
\omega(t,x)=\sum_k\,a_k(0)\,e^{(\lambda_k + \delta)\,t}\,\psi_k(x).
\]
Since $\alpha(s)$ is irreducible, the lift $u$ such that $u^*\alpha(s)
= \alpha(s)$ is determined uniquely up to a sign, so that the form
$\omega$ satisfies the periodic boundary condition $\omega (t+1,x) = u^*\,
\omega (t,x)$. On every eigenspace of $K_{\alpha(s)}$ with a fixed
eigenvalue $\lambda$, this translates into
\[
u^*\,\sum_{\lambda_k = \lambda}\,a_k(0)\,\psi_k(x) =
e^{\lambda + \delta}\cdot\sum_{\lambda_k = \lambda}\,a_k(0)\,\psi_k(x).
\]
This means that $e^{\lambda + \delta}$ is an eigenvalue of $u^*$. Since
$u^*$ has finite order and $\lambda + \delta$ is a real number, we must
have $\lambda + \delta = 0$.  Therefore, $\omega$ belongs to the kernel
of $K^u_{\alpha(s)} + \delta$. Conversely, it is easy to see that $\ker
(D_{A(s)} - \delta)$ is contained in $\ker (K^u_{\alpha(s)} + \delta)$.

A similar argument using the formula $D^*_{A(s)} = -\p/\p t - K_{\alpha(s)}$
shows that $\coker (D_{A(s)} - \delta) = \ker (K^u_{\alpha(s)} + \delta)$.

For a small generic $\delta > 0$, the spectral curves of $K^u_{\alpha(s)}
+ \delta$ intersect the $s$--axis transversally. Since $H_s = D_{A(s)} -
D_{\theta} = K^u_{\theta} - K^u_{\alpha(s)} = - h_s$,  this implies that
the homotopy $H_s$ is regular. The proof is now complete.
\end{proof}

\begin{remark}\label{R:hom}
A similar argument can be used to show that $H^1_{\tau}(\Sigma,\;\ad\alpha)
= H^1 (X;\ad A)$ thus providing an independent proof of Proposition \ref{P2}.
\end{remark}

To finish the orientation comparison, we need to evaluate $\ep(D_{\theta},
H, D_A)$.  To this end, we choose a path of operators from $D_{\theta}$ to
$D_A$ which consists of three segments. The first segment is $D_{\theta}
- s\cdot\delta$, $0\le s\le 1$, connecting $D_{\theta}$ to $D_{\theta} -
\delta$, the second is $D_{A(s)} - \delta$ as in Proposition \ref{P:delta},
and the last is $D_A - (1-s)\cdot \delta$, $0\le s\le 1$, connecting $D_A
- \delta$ to $D_A$. Since the orientation transport is additive, the
orientation transports $\ep(D_{\theta},H, D_A)$ and $\ep(D_{\theta} -
\delta, H, D_A - \delta)$, for a small generic $\delta > 0$, differ by the
product of orientation transports along the first and the last segments.
If $\delta$ is small enough, the operators in the last segment have zero
kernels, making for the trivial orientation transport. The family
$D_{\theta} - s\cdot\delta$ only has kernel at $s = 0$.  That kernel is
isomorphic to $H^0 (X;\ad\theta) = \su(2)$, with the resonance operator
$R_0 = -\delta\cdot\Id\co  \su(2) \to \su(2)$ so that $\det R_0 = -\delta^3
< 0$.  Finally, note that in formula~\eqref{eq} for the orientation
transport, the product of the terms $(-1)^{d_s}$ starts at $s=0$.  Since
$\dim\ker D_{\theta} = 3$, we conclude that the orientation transport
along the first segment is also equal to $+1$, and $i^*$ is orientation
preserving.


\subsection{A note on the Lefschetz number}\label{S:lef}
Under certain non-degeneracy assumptions on the character variety $\R^*
(\Sigma)$, one can express the Furuta--Ohta invariant of a mapping torus
in terms of the action of $\tau$ on the Floer homology of $\Sigma$. The
standing assumption we make in this subsection is as follows (compare
with Section~\ref{S:perturb})\,:
\begin{enumerate}
\item[{$(*)$}] There is an equivariant admissible perturbation $h$ such
that the perturbed moduli space $\R_h^* (\Sigma)$ is non-degenerate.
\end{enumerate}
Note that this condition is not as strong as requiring that
$\R^*(\Sigma)$ be non-degenerate, yet it is stronger than is guaranteed
by the construction in Section~\ref{S:perturb}. We will see in
Section~\ref{S:examples} that condition $(*)$ holds in some interesting
examples.

Under the above assumption, there is an action of $\tau^*$ on $\R_h^*
(\Sigma)$.  For an automorphism $f\co  V \to V$ of a graded module $V$, denote
by $L(f,V)$ its Lefschetz number. Note that this only requires a
$\Z_2$--grading, and so makes sense for the action of $\tau$ on the
$\Z_8$--graded Floer chain complex $IC_*(\Sigma)$ and Floer homology groups
$I_*(\Sigma)$.
\begin{theorem}\label{T:lefschetz} Let $\tau\co  \Sigma \to \Sigma$ be an
orientation-preserving diffeomorphism such that the non-degeneracy
condition $(*)$ holds.  Then the Furuta--Ohta invariant of the mapping torus
$X$ of $\tau$ is given by
$$
\lambda_{FO}(X)= \frac12\, L(\tau,I_*(\Sigma)).
$$
\end{theorem}
\begin{proof}
Write $\R^*_h(\Sigma) = \R^p \cup \R^\tau$ where $\R^\tau$
is the fixed point set of $\tau^*$ acting on $\R^*_h(\Sigma)$, and
points in $\R^p$ are permuted by $\tau^*$. The submodules of the Floer
chain complex $IC_*(\Sigma)$ generated by these sets will be called
$IC_*^p$ and $IC_*^{\tau}$. They need not be sub chain complexes;
however, the chain map $\tau_*\co  IC_*(\Sigma) \to IC_*(\Sigma)$ induced
by $\tau_*$ preserves these two submodules, as well as the grading.

It follows that (the first equality is the Hopf trace formula; the
second is just linear algebra)
\begin{equation}\label{E:hopf}
L(\tau, I_*(\Sigma)) = L(\tau,IC_*(\Sigma)) =
L(\tau,IC_*^p) + L(\tau,IC_*^{\tau}).
\end{equation}
Note that by construction, the diagonal entries for $(\tau,IC_*^p)$
are all $0$, and the ones for $(\tau,IC_*^{\tau})$ are $\pm 1$, whose
signs we now have to figure out.

Let us follow \cite[Section 5.4]{donaldson:floer} and orient every
generator $\alpha$ in the Floer chain complex by the respective line
bundle $\lambda_W (\alpha)\,\otimes\,\lambda^*_W (\theta)$, where
$W$ is a smooth compact oriented 4--manifold with boundary $\Sigma$,
and $\lambda_W (\alpha)$ is the determinant bundle over the moduli
space of connections on $W$ with limiting value $\alpha$. Now, to
compute the sign of $\tau^* \alpha$, we simply attach the mapping
cylinder $C$ of $\tau$ to $W$. Then $\tau^* \alpha$ is oriented by
$\lambda_{W\cup C}(\tau^*\alpha)\,\otimes\,\lambda_{W\cup C}(\theta)$.
By the excision principle, we have
\[
\lambda_{W\cup C} (\tau^*\alpha) = \lambda_W (\alpha)\,\otimes\,
\lambda_C (\alpha,\tau^*\alpha)
\]
and $\lambda_{W\cup C} (\theta) = \lambda_W (\theta)\,\otimes\,1$.
Thus the sign of the diagonal entry in $(\tau,IC_*^{\tau})$
corresponding to $\alpha$ is $\lambda_C (\alpha,\tau^*\alpha)$,
which, by excision principle, also orients the respective point in
$\M^*(X)$ (compare with Section \ref{S:orient4}).

Taking into account the two-to-one correspondence between $\M^*(X)$
and $\R^{\tau}(\Sigma)$ (given by Proposition \ref{P1} and, in the
perturbed case, by Proposition \ref{P:pert21}), we obtain the result.
\end{proof}

\begin{remark} The observation that a gauge-theoretic invariant of a
mapping torus may be interpreted as a Lefschetz number occurs in several
other papers~\cite{frohman-nicas,boden-nicas,salamon:lefschetz}. These
works are concerned with mapping tori of surface diffeomorphisms, rather
than $4$--dimensional mapping tori.  Because the moduli spaces of unitary
connections (resp. solutions to the Seiberg--Witten equations) on a surface
are well-understood, the non-degeneracy assumptions in~\cite{frohman-nicas,
boden-nicas} (resp.~\cite{salamon:lefschetz}) are milder than our
condition $(*)$.

We conjecture that the equality of $\lambda_{FO}$ and the Lefschetz number
holds without any additional smoothness assumptions.  It would be of some
interest, with regard to the conjecture of Furuta and Ohta, to see some
direct relationship between the Lefschetz number and the Rohlin invariant.
\end{remark}


\section{Definition of the Furuta--Ohta invariant}\label{D:fo}
Before we go on to prove Theorem \ref{T:finite} in general, we give a
detailed description of the class of perturbations needed to properly define
the Furuta--Ohta invariant, as such a description is not detailed in
\cite{furuta-ohta}. The proof of Theorem \ref{T:finite} in the next section
will then amount to finding sufficiently many \emph{equivariant}
perturbations in this class.


\subsection{Admissible perturbations}
Let $X$ be a $\Z[\Z]$--homology $S^1\times S^3$. Consider an embedding
$\psi\co  S^1 \to X$ and extend it to an embedding  $\psi\co  S^1\times N^3
\to X$ where $N^3$ is an oriented 3--manifold.  Thus for each
point $x\in N^3$ we have a parallel copy $\psi(S^1\times\{x\})$ of
$\psi\co  S^1\to X$.

Let $\bar P$ be a trivial $SU(2)$--bundle over $X$, and let $\A(X)$ be the
affine space of connections in $\bar P$. Given $A\in \A(X)$, denote by
$\hol_A (\psi (S^1\times\{x\}),s)\in SU(2)$ the holonomy of $A$ around
the loop $\psi (S^1\times\{x\})$ starting at the point $\psi (s,x)$.
Let $\Pi\co  SU(2)\to \su(2)$ be the projection given by
\[
\Pi (u) = u - \frac 1 2\,\tr(u)\cdot\Id.
\]
It is equivariant with respect to the adjoint action on both $SU(2)$
and $\su(2)$, therefore, assigning $\Pi\hol_A(\psi(S^1\times\{x\}),s)$
to $\psi(s,x)\in X$ defines a section of $\ad\bar P = \bar P
\times_{\ad}\su(2)$ over $\psi (S^1\times N^3)$. Let $\nu\in\Omega^2_+
(X)$ be a self-dual 2--form supported in $\psi (S^1\times N^3)$ and
define a section
\[
\sigma (\nu, \psi, A) \in \Omega^2_+ (X,\ad\bar P)
\]
by taking tensor product of $\Pi \hol_A (\psi (S^1\times\{x\}),s)$ with
$\nu$ over $\psi (S^1\times N^3)$ and letting it be zero otherwise. For
fixed $\nu$ and $\psi$, this defines a map $\sigma\co  \A(X) \to \Omega^2_+
(X, \ad\bar P)$ which is equivariant with respect to the gauge group
$\G(X)$.

More generally, let $\Psi = \{\,\psi_k\,\}$ be a collection of embedded
loops $\psi_k\co  S^1\to X$, $k = 1,\ldots, n$, with disjoint images. We
will refer to $\Psi$ as a \emph{link}. Extend $\Psi$ to a collection of
embeddings $\psi_k\co  S^1\times N^3_k\to X$ as above so that the
$\psi_k (S^1\times N^3_k)$ are disjoint. For any choice of $n$ smooth
functions $\bar f_1,\ldots,\bar f_n\co  [-2,2]\to\mathbb R$ with vanishing
derivatives at $\pm\, 2$, define $\sigma\co  \A(X) \to \Omega^2_+ (X, \ad
\bar P)$ by the formula
\begin{equation}\label{E:sigma}
\sigma (A) = \sum_{k=1}^n\; \p \bar f_k\cdot\sigma(\nu_k, \psi_k, A),
\end{equation}
where $\p \bar f_k$ is the function $\bar f'_k$ evaluated at $\tr
\hol_A (\psi_k (S^1\times\{x\}),s)$, and $\nu_k$ are real valued
self-dual forms on $X$, each supported in its respective $\psi_k
(S^1\times N^3_k)$. We call $\sigma$ an \emph{admissible perturbation
relative to $\Psi$}.

For a fixed choice of the forms $\nu_1,\ldots,\nu_n$, denote by
$\F_{\Psi}$ the space of admissible perturbations relative to $\Psi$
with the $C^3$--topology given by the correspondence $\sigma \mapsto
(\bar f_1, \ldots, \bar f_n)$.


\subsection{Perturbed ASD connections}
Let $\sigma\co  \A(X) \to \Omega^2_+ (X, \ad\bar P)$ be an admissible
perturbation and define $\zeta_{\sigma}\co  \A(X)\to\Omega^2_+(X,\ad\bar P)$
by the formula
\[
\zeta_{\sigma} (A) = F_+ (A) + \sigma (A).
\]
A connection $A$ is called \emph{perturbed ASD} if $\zeta_{\sigma}(A)
= 0$.  The moduli space of perturbed ASD connections will be denoted by
$\M_{\sigma}(X)$ so that $\M_{\sigma}(X) = \zeta_{\sigma}^{-1}(0)/\G(X)$.

If $\sigma = 0$ then $\M_{\sigma}(X)$ coincides with the flat moduli
space $\M(X)$. Since $\pi_1 (X)$ is finitely presented, $\M(X)$ is a
compact real algebraic variety. According to \cite{furuta-ohta}, reducible
flat connections form a connected component in $\M(X)$, therefore the moduli
space $\M^*(X)$ of irreducible flat connections is also compact. We are
only interested in irreducible connections hence we will not perturb
the reducible part of $\M(X)$. Note that small enough perturbations do
not create new reducible connections. Denote by $\M^*_{\sigma}(X)$ the
moduli space of irreducible perturbed ASD connections; for a sufficiently
small $\sigma$, it is compact.

The local structure of $\M^*_{\sigma}(X)$ near a point $[A] \in
\M^*_{\sigma}(X)$ is described by the deformation complex
\[
\begin{CD}
\Omega^0 (X, \ad\bar P) @> d_A >> \Omega^1 (X, \ad\bar P)
@> d^+_A + D\sigma (A) >> \Omega^2_+ (X, \ad\bar P). \\
\end{CD}
\]

We call $\M^*_{\sigma}(X)$ \emph{non-degenerate} at $[A]\in \M^*_{\sigma}
(X)$ if the second cohomology of the above complex vanishes, that is,
$H^2_{\sigma,+}(X;\ad A) = \coker(d_A^+ + D\sigma(A)) = 0$. Since $A$
is irreducible, an equivalent condition is vanishing of
$\coker(D_A^{\,\sigma})$, where $D_A^{\,\sigma} = D_A + D\sigma (A)$ is
the perturbed ASD operator (\ref{E:da}). We call $\M^*_{\sigma}(X)$
\emph{non-degenerate} if it is non-degenerate at all $[A]\in\M^*_{\sigma}
(X)$. If $\M^*_{\sigma}(X)$ is non-degenerate, it consists of finitely
many points oriented by the orientation transport $\ep (D_{\theta}^{\,\sigma},
H,D_A^{\,\sigma}) = \ep (D_{\theta},H,D_A^{\,\sigma})$, compare with
Section \ref{S:orient4}.


\subsection{Abundance of perturbations}
We wish to show that we have enough admissible perturbations $\sigma$ to
find a non-degenerate $\M^*_{\sigma}(X)$. A link $\Psi = \{\,\psi_k\,\}$
will be called \emph{abundant} at an ASD connection $A$ if there exist
self dual 2--forms $\nu_k$ supported in $\psi_k (S^1\times N^3_k)$ such
that the sections $\sigma(\nu_k,\psi_k, A)$ span the vector space $H^2_+
(X;\ad A)$. Note that if a link $\Psi$ is abundant at $A$ it is also
abundant at any ASD connection gauge equivalent to $A$.

\begin{lemma}\label{L:abundance}
There exists a link $\Psi$ which is abundant at all $[A]\in\M^*(X)$.
\end{lemma}

\begin{proof}
That there is a link which is abundant at any given $[A] \in \M^*(X)$ is
immediate from Lemma (2.5) of \cite{donaldson:orientation}. Abundance is
an open condition with respect to $A$ hence existence of the universal
link $\Psi$ follows from compactness of $\M^*(X)$.
\end{proof}

\begin{proposition}
For a small generic admissible perturbation $\sigma$, the moduli space
$\M_{\sigma}^* (X)$ is non-degenerate.
\end{proposition}

\begin{proof}
Let us fix an abundant link $\Psi$ and consider the map $\Phi\co  \B^* (X)
\times \F_{\Psi} \to \Omega^2_+ (X, \ad\bar P)$ given by the formula
\[
\Phi(A, \bar f) = F_+(A) + \sum_{k=1}^n\;\p \bar f_k \cdot \sigma(\nu_k,
\psi_k,A)
\]
where $\bar f = (\bar f_1,\ldots, \bar f_n) \in \F_{\Psi}$. The partial
derivative of this map with respect to $A$ has cokernel $H^2_+(X;\ad A)$,
while the partial derivative with respect to $\bar f$ is onto $H^2_+
(X;\ad A)$ according to Lemma \ref{L:abundance}.

Therefore, $\Phi$ is a submersion; in particular, $\Phi^{-1}(0)$ is
locally a smooth manifold. Consider projection of $\Phi^{-1} (0)\subset
\B^*(X)\times \F_{\Psi}$ onto the second factor. By the Sard--Smale
theorem, this projection is a submersion for a dense set of $\bar f$
within a sufficiently small neighborhood of zero in $\F_{\Psi}$. Since
the pre-image of $\bar f$ is the moduli space $\M_{\sigma}^*(X)$ with
$\sigma (A) = \sum\;\p \bar f_k \cdot \sigma(\nu_k,\psi_k,A)$, we are
finished.
\end{proof}


\subsection{Definition of the invariant}
The Furuta--Ohta invariant is defined as one fourth of the signed count of
points in a non-degenerate moduli space $\M_{\sigma}^*(X)$. That it is
well defined follows from the standard cobordism argument\,: for any
generic path of small admissible perturbations $\sigma(t)$, $t\in [0,1]$,
the parameterized moduli space
\[
\overline{\M}^{\,*} (X) = \bigcup_{t\in [0,1]}\; \M^*_{\sigma(t)}(X)
\times \{\,t\,\}
\]
is a \emph{compact} oriented cobordism between $\M^*_{\sigma(0)}(X)$ and
$\M^*_{\sigma(1)}(X)$. The cobordism $\overline{\M}^{\,*}(X)$ is compact
because the condition $H_* (\tilde X; \Z) = H_* (S^3; \Z)$ ensures that
the reducible connections are isolated in the flat moduli space, see
\cite{furuta-ohta}.


\section{The perturbed case}\label{S:perturb}
In this section, we prove Theorem \ref{T:finite} without assuming that
$\R^{\tau}(\Sigma)$ is non-degenerate. We begin by recalling the
definition of the equivariant Casson invariant, which was given in
\cite{collin-saveliev:casson} for finite order diffeomorphisms $\tau$
having fixed points, and then extending it to cover the fixed point free
$\tau$ (still of finite order). Theorem \ref{T:finite} follows after
matching perturbations in the two theories.


\subsection{Perturbing $\R^{\tau}(\Sigma)$}
The decomposition (\ref{E:decomp2}) splits $\R^{\tau} (\Sigma)$ into a
union of $\R^u (\Sigma) = \R^{\tau}(\Sigma) \cap \B^u (\Sigma)$ each of
which consists of (the gauge equivalence classes of) flat connections
$\alpha$ such that $u^*\alpha = \alpha$, for a constant lift $u\co  P \to
P$. Clearly, it is sufficient to do perturbations on each of the
$\B^u (\Sigma)$ separately, hence we will fix a constant lift $u$ from
the beginning.

Let $\{\,\gamma_j\,\}$ be a link in $\Sigma$, that is, a collection
$\gamma_j\co  S^1\times D^2\to\Sigma$, $j = 1,\ldots, n$, of embeddings
with disjoint images. For any $\alpha \in \A (\Sigma)$, denote by
$\hol_{\alpha}(\gamma_j (S^1\times\{z\}),s) \in SU(2)$ the holonomy
of $\alpha$ around the loop $\gamma_j (S^1\times \{z\})$ starting at
the point $\gamma_j (s,z)$. Let $\eta(z)$ be a smooth rotationally
symmetric bump function on the disc $D^2$ with support away from the
boundary of $D^2$ and with integral one, and let $f_j\co  SU(2) \to
\mathbb R$ be a collection of smooth functions invariant with respect
to conjugation. Following \cite{taubes:casson}, define an
\emph{admissible perturbation} $h\co  \A^u(\Sigma)\to\mathbb R$ by the
formula
\begin{equation}\label{ha}
h(\alpha) = \sum_{j=1}^n\;\int_{D^2}\;f_j(\hol_{\alpha}(\gamma_j (S^1
\times\{z\}),s))\,\eta(z)\,d^2 z.
\end{equation}
The conjugation invariance of the $f_j$ ensures that the function $h$
does not depend on $s$. The action of the gauge group only changes
holonomies around $\gamma_j (S^1\times\{z\})$ within their $SU(2)$
conjugacy classes. Therefore, the admissible perturbation $h$ is well
defined on $\B^u (\Sigma)$.

Suppose that the link $\{\,\gamma_j\,\}$ is invariant with respect to the
action of $\tau$ and that $f_j = f_m$ whenever $\tau_*\gamma_j=\gamma_m$.
Then formula (\ref{ha}) defines an admissible perturbation $h$ such that
$h(\tilde\tau^*\alpha)= h(\alpha)$ for any lift $\tilde\tau$. We will call
$h$ an \emph{equivariant admissible perturbation}. The formula
\[
\zeta^u_h (\alpha) = *\,F (\alpha) + \nabla h (\alpha),
\]
where $\nabla h$ is the $L^2$--gradient of $h$, defines a map $\zeta^u_h\co 
\A^u (\Sigma) \to \Omega^1_u (\Sigma,\ad P)$ where $\omega \in \Omega^1_u
(\Sigma,\ad P)$ if and only if $u^*\omega =\omega$, see Section \ref{S:R}.

A connection $\alpha \in \A^u(\Sigma)$ is called \emph{perturbed flat} if
$\zeta^u_h (\alpha) = 0$. The moduli space of the perturbed equivariant
flat connections will be denoted by $\R^u_h (\Sigma)$ so that $\R^u_h
(\Sigma) = (\zeta_h^u)^{-1}(0)/\G^u(\Sigma)$. If $h = 0$ then $\R^u_h
(\Sigma)$ coincides with the moduli space $\R^u (\Sigma)$ of equivariant
flat connections.

The moduli space $\R^u_h(\Sigma)$ is called \emph{non--degenerate} if
$H^1_{h,\tau} (\Sigma;\ad\alpha) = \ker (d_{\alpha}\zeta^u_h)\allowbreak
\cap\ker d^*_{\alpha}$ vanishes for all $[\alpha] \in \R^u_h (\Sigma)$.
After a minor modification to cover the case of empty $\Fix(\tau)$,
the argument of \cite{collin-saveliev:casson}, Section 3.8, shows that
there exist small equivariant admissible perturbations $h\co \A(\Sigma)\to
\mathbb R$ making all the moduli spaces $\R^u_h(\Sigma)$ non-degenerate.
For such a perturbation, the space $\R^{\tau}_h(\Sigma)$, which is a
union of the $\R^u_h(\Sigma)$, is smooth compact manifold of dimension
zero.


\subsection{Definition of the equivariant Casson invariant}\label{D:eq}
Let $\R^{\tau}(\Sigma)$ be the equivariant flat moduli space, and choose
(if necessary) a small equivariant admissible perturbation $h$ such that
$\R^{\tau}_h (\Sigma)$ is non-degenerate. Define the \emph{equivariant
Casson invariant} by the formula
\[
\lambda^{\tau}(\Sigma) = \frac 1 2\;\sum_{\alpha\in\R^{\tau}_h(\Sigma)}
(-1)^{\mu^{\tau}(\alpha)},
\]
where $\mu^{\tau}(\alpha)$ is the equivariant spectral flow of the perturbed
odd signature operators $K_{h,\alpha(s)} = K_{\alpha(s)}+\Hess_{\,\alpha(s)}
h$ as in Section \ref{S:R}. If $\tau$ has fixed points, this invariant
coincides with the equivariant Casson invariant defined in
\cite{collin-saveliev:casson}.

\begin{proposition}
The invariant $\lambda^{\tau}(\Sigma)$ is well defined as an integer
valued invariant of the pair $(\Sigma, \tau)$.
\end{proposition}

\begin{proof}
This was proved in \cite{collin-saveliev:casson} in the case when $\tau$
has fixed points. For a fixed point free $\tau$, the result is immediate
from Theorem \ref{thm2}, which expresses $\lambda^{\tau}(\Sigma)$ in terms
of the Boyer--Nicas type invariants of $\Sigma/\tau$.
\end{proof}


\subsection{Matching the perturbations}
Let $X$ be the mapping torus of $\tau\co \Sigma \to \Sigma$. If $\M^* (X)$
happens to be degenerate then so is the flat moduli space $\R^\tau(\Sigma)$,
see Proposition \ref{P2}. Our goal in this section will be to perturb
$\M^* (X)$ and $\R^\tau(\Sigma)$ in a consistent manner so as to keep the
two-to-one correspondence between them which existed in the unperturbed
case.

Let us fix a constant lift $u\co  P \to P$ and let $\bar P$ be the trivial
$SU(2)$--bundle over $X$ obtained by pulling $P$ back to $[0,1]\times
\Sigma$ and identifying the ends via $u$. The maps $\pi\co  \A^u (\Sigma)
\to \A^* (X)$, see Remark \ref{R:pi}, can be included into the
commutative diagram
\[
\begin{CD}
\A^u (\Sigma) @> *\,F(\alpha) >> \Omega^1_u (\Sigma,\ad P) \\
@VV \pi V  @VV p V \\
\A^* (X) @> F_+ (\pi(\alpha)) >> \Omega^2_+ (X, \ad\bar P) \\
\end{CD}
\]

where $p$ is defined as follows. Given $\omega\in\Omega^1_u(\Sigma,\ad P)$,
pull the 2-form $*\,\omega$ back to $[0,1]\times\Sigma$ and identify the
ends via $u^*$. The form $p(w)$ is then the self--dual part of this 2-form.
The zero set of $*\,F(\alpha)$ is the flat moduli space $\R^u(\Sigma)$, and
$\pi$ establishes a one-to-two correspondence between the union of all the
$\R^u(\Sigma)$ and the moduli space $\M^* (X)$, see Proposition \ref{P1}.

\begin{proposition}\label{P:4.1}
There exists an equivariant admissible perturbation $h$ and an admissible
perturbation $\sigma$ such that, for every $u$, the diagram
\begin{equation}\label{E:cd}
\begin{CD}
\A^u (\Sigma) @> \zeta^u_h(\alpha) >> \Omega^1_u(\Sigma,\ad P) \\
@VV \pi V  @VV p V \\
\A^* (X) @> \zeta_{\sigma}(\pi(\alpha)) >> \Omega^2_+(X,\ad\bar P) \\
\end{CD}
\end{equation}

commutes, and both moduli spaces $\R^{\tau}_h(\Sigma)$ and $\M_{\sigma}^*
(X)$ are non-degenerate.
\end{proposition}

\begin{proof}
We start with an equivariant admissible perturbation $h\co  \A^u(\Sigma)
\to {\mathbb R}$ and evaluate explicitly $\zeta^u_h
(\alpha) = * \,F(\alpha) + \nabla h (\alpha)$. Write each of the $f_j\co 
SU(2)\to \mathbb R$ as $f_j = \bar f_j\circ\tr$ for a smooth function
$\bar f_j \co  [-2,2]\to \mathbb R$.  Then, according to
\cite{herald:perturbations}, the gradient $\nabla h(\alpha)$ is the
1-form
\[
\nabla h (\alpha) = \sum_{j = 1}^n\;
\bar f'_j (\tr\hol_{\alpha}(\gamma_j (S^1\times\{z\},s)))\cdot
\Pi\hol_{\alpha} (\gamma_j (S^1 \times\{z\}),s) \,\eta(z) \,ds
\]
where each of the summands is understood to have support in its respective
copy of $\gamma_j (S^1\times D^2)$. Next, each orbit of the $\Z_n$--action
induced by $\tau$ on $\{\,\gamma_j\,\}$ gives rise to a copy of $S^1\times
S^1\times D^2$ and hence to the link $\psi_k\co  S^1 \times N^3_k \to X$ with
$S^1$--direction corresponding to $\gamma_j (S^1\times\{0\})$ and $N^3_k =
S^1\times D^2$. For any connection $A$ over $X$ set $\nu_k=p(\eta(z)\,ds)$
and define $\sigma$ as in (\ref{E:sigma}).  An easy calculation shows that
$\hol_{\pi(\alpha)}(\psi_k (S^1\times\{(t,z)\}), s) = \hol_{\alpha}
(\gamma_j (S^1\times \{z\}),s)$, for a respective $j$ in the $\Z_n$--orbit
which gives rise to $\psi_k$. Therefore, $\sigma$ makes the diagram
(\ref{E:cd}) commute.

The second statement follows from the fact that the map $p^*\co  H^1_{h,\tau}
(\Sigma; \ad\alpha) \to H^2_{\sigma,+} (X; \ad\pi(\alpha))$ is an
isomorphism. This follows via a minor modification of the proof of
Proposition \ref{P:delta}, the crucial observation being that
$D^{\sigma}_A = \p/\p t - K_{h,\alpha}$, compare with Remark \ref{R:hom}.
\end{proof}

Let us fix $h$ and $\sigma$ as in Proposition \ref{P:4.1}. The following
result is a perturbed analogue of Proposition \ref{P1}.

\begin{proposition}\label{P:pert21}
The map $\pi$ provides a one-to-two correspondence between $\R^{\tau}_h
(\Sigma)$ and $\M^*_{\sigma}(X)$.
\end{proposition}

\begin{proof}
We only need to check that all connections in $\M^*_{\sigma}(X)$ are of
the form $\pi(\alpha)$, up to gauge equivalence. Given a connection $A$
over $X$, put it into temporal gauge over $[0,1] \times \Sigma$ so that
$A = \{\,A(t)\,\}$. The connections $A(0)$ and $A(1)$ will be gauge
equivalent (via the holonomy along $[0,1]\times\{\,x\,\}$) although not
necessarily equal. Let $\delta (\alpha)$ be the pull-back of $*\,\nabla
h (\alpha)$ over each of the product regions $[0,1]\times \gamma_j (S^1
\times D^2)$, and zero otherwise.  Then $\sigma(A) = \delta (A)_+$ so
that the connection $A$ is perturbed ASD if and only if $F_+ (A) +
\delta_+ (A) = 0$.

We wish to show first that $A$ is perturbed flat, that is, $F (A) +
\delta (A) = 0$. To this end, calculate

\begin{alignat}{1}
\| F(A) + \delta (A)\|^2
&= -\int_X \tr (F(A)+\delta(A))\wedge
*\, (F(A)+\delta(A))  \notag \\
&= \quad\int_X \tr (F(A)+\delta(A))\wedge
(F(A)+\delta(A)) \notag \\
&= \quad\int_X \tr (F(A)\wedge F(A)) +
2\,\int_X \tr (F(A)\wedge\delta(A)) \notag
\end{alignat}
since $\delta(A)\wedge\delta(A) = 0$. The first integral vanishes by the
Chern--Weil theory. To compute the second, pull back to $[0,1] \times
\Sigma$ where we write $F(A) = dt \wedge A'(t) + F(A(t))$ so that $F(A)
\wedge\delta(A) = dt\wedge A'(t)\wedge\delta(A)$. First we integrate over
$\Sigma$ to obtain
\[
\int_{\Sigma} \tr (F(A)\wedge\delta(A)) =
\int_{\Sigma} \tr (A'(t)\wedge *\, \nabla h (A(t))) =
\frac d {dt}\, h(A(t))\;dt.
\]
Integration with respect to $t$ now yields
\[
\int_X \tr (F(A)\wedge\delta(A)) =
\int_0^1 \frac d {dt}\,h(A(t))\,dt
= h(A(1)) - h(A(0)) = 0
\]
since $h$ is invariant with respect to gauge transformations.

Now, since $A$ is perturbed flat, it satisfies the equation
$dt\wedge A'(t) + F(A(t)) + *\,\nabla h (A(t)) = 0$ over
$[0,1]\times\Sigma$ with periodic boundary conditions. Therefore,
$A(t) = \alpha$ for all $t$, where $\alpha \in \A^u (\Sigma)$ is
a perturbed flat connection, so that $A = \pi(\alpha)$ after a
gauge transformation if necessary.
\end{proof}


\subsection{Proof of Theorem \ref{T:finite}}
Choose perturbations $h$ and $\sigma$ as in Proposition \ref{P:4.1},
then $\M_{\sigma}^*(X)$ and $\R^{\tau}_h (\Sigma)$ are in two-to-one
correspondence and are both non-degenerate. Moreover, in temporal
gauge,  $D^{\sigma}_A  =  \p/\p t - K_{h,A(t)}$ where $K_{h,A(t)}$
is the perturbed operator $K_{A(t)}$. The same argument as in the
unperturbed case now shows that that the corresponding points in
$\M^*_{\sigma} (X)$ and $\R^{\tau}_h (\Sigma)$ are counted with the
same sign. This completes the proof of Theorem \ref{T:finite}.


\section{The equivariant Casson invariant:\; $\Fix(\tau) \ne \emptyset$}
\label{S6}
In the case of non-empty $\Fix(\tau)$, the equivariant Casson invariant
can be easily evaluated using an explicit formula proved in
\cite{collin-saveliev:casson}. In
this case, the quotient space $\Sigma' = \Sigma/\tau$ is an integral
homology sphere, and the projection $\Sigma\to \Sigma'$ is an $n$--fold
branched cover with branch set a knot $k \subset \Sigma'$. According to
Theorem 1 of \cite{collin-saveliev:casson},
\begin{equation}\label{E:cs}
\lambda^{\tau}(\Sigma) = n\cdot\lambda(\Sigma') + \frac 1 8\,
\sum_{m = 0}^{n-1}\; \sign^{m/n} k,
\end{equation}
where $\sign^{\alpha} k$ is the \emph{Tristram-Levine equivariant knot
signature}, defined as the signature of the Hermitian form
$(1 - \exp(2\pi i\alpha))\,S + (1 - \exp(-2\pi i\alpha))\,S^t$ for any
choice of a Seifert matrix $S$ of $k$.


\section{The equivariant Casson invariant:\; $\Fix(\tau) = \emptyset$}
In this section, we obtain a formula for $\lambda^{\tau}(\Sigma)$ in the
case of a fixed point free $\tau$. We first express $\lambda^{\tau}(\Sigma)$
in terms of certain Boyer-Nicas type invariants (which extend invariants of
\cite{boyer-nicas} and \cite{boyer-lines}), and then identify the latter
via topology. We will be using interchangeably $G$--gauge theories with $G =
SO(3)$, $SU(2)$, and $U(2)$; to keep better track of $G$ we will include it
in our notations -- for example, $\R(\Sigma,SU(2))$ will mean the moduli
space of flat $SU(2)$--connections etc.


\subsection{Identifying moduli spaces}
Let $P$ be a trivialized $SU(2)$--bundle over $\Sigma$.  According to
Section \ref{S:decomp3}, the equivalence classes of lifts $\tilde\tau\co
P\to P$ are classified by the group $H^2 (\Sigma';\Z_2)$. The
equivalence class of lifts corresponding to $w\in H^2(\Sigma';\Z_2)$
will be denoted by $\tau_w$. The decomposition (\ref{E:decomp1}) then
takes the form
\[
\B^{\tau}(\Sigma,SU(2)) = \bigsqcup_w\; \B^{\tau}_w (\Sigma,SU(2)),
\]
where $\B^{\tau}_w(\Sigma,SU(2))$ consists of irreducible
$SU(2)$--connections $A$ in $P$ such that $\tau^*_w A = A$, modulo the
group of gauge transformations $g\co  P \to P$ such that $g\tau_w = \pm\,
\tau_w g$. Since $H^2(\Sigma';\Z_2)$ is trivial if $n$ is odd and is
equal to $\Z_2$ if $n$ is even, we have at most two different equivalence
classes of lifts, $\tau_0$ and $\tau_1$. We will fix constant lifts
within each class so that $\tau_0 (x,y) = (\tau(x),y)$ and $\tau_1(x,y) =
(\tau(x),u\cdot y)$ where $u\in SU(2)$ is such that $u^n = -1$, see
Section \ref{S:decomp3}.

For any choice of $\tau_w$, the quotient bundle $P'_w = P/\tau_w$ over
$\Sigma' = \Sigma/\tau$ is an $SO(3)$--bundle with $w_2(P'_w) = w$. It
has a natural smooth structure such that the quotient map $\pi\co  P \to
P'_w$ is a smooth bundle morphism. Denote by $\B^*_w(\Sigma',SO(3))$
the space of the gauge equivalence classes of irreducible connections
in $P'_w$. Let $\R^*_w(\Sigma',SO(3))\subset\B^*_w(\Sigma',SO(3))$ and
$\R^{\tau}_w (\Sigma,SU(2))\subset \B^{\tau}_w(\Sigma,SU(2))$ be the
respective moduli spaces of irreducible flat connections.

\begin{proposition}\label{P:pullback}
The pull back of connections via $\pi$ induces a homeomorphism between
$\R^*_w (\Sigma',SO(3))$ and $\R^{\tau}_w(\Sigma,SU(2))$.
\end{proposition}

\begin{proof}
First note that the $SU(2)$ and $SO(3)$ gauge theories on $\Sigma$ are
equivalent because $\Sigma$ is an integral homology sphere. The result
then essentially follows from the discussion in Section \ref{S:decomp3},
after we check that the pull back $\pi^*\alpha'$ of an irreducible flat
connection $\alpha'$ in $P'_w$ is irreducible. Let $\alpha'\co  \pi_1
(\Sigma') \to SO(3)$ be the holonomy representation of $\alpha'$. Its
pull back via $\pi_1(\Sigma)\to \pi_1(\Sigma')$ is the holonomy
representation of $\pi^*\alpha'$. Suppose that $\pi^* \alpha'$ is
reducible, then $\pi^*\alpha'$ is trivial because $\Sigma$ is an
integral homology sphere. This means that $\alpha'$ factors through
$\pi_1(\Sigma')/\pi_1(\Sigma) = \Z_n$ and hence is also reducible, a
contradiction.
\end{proof}

\begin{remark}
If $n$ is odd, we have just one $SO(3)$--bundle on $\Sigma'$ which is
trivial and hence admits a unique $SU(2)$--lift. In this case,
Proposition \ref{P:pullback} also establishes a bijective correspondence
at the level of flat $SU(2)$--connections.
\end{remark}


\subsection{The Boyer-Nicas type invariants}
Let $\Sigma' = \Sigma/\tau$ as before.
Evidently, $\Sigma'$ is \emph{cyclically finite} in the terminology
of~\cite{boyer-nicas}, ie, every cyclic cover of $\Sigma'$ is a rational
homology sphere. In this situation, Boyer and Nicas defined an integer
valued invariant, by counting irreducible $SU(2)$--representations, with
signs defined using a Heegaard splitting as in Casson's original
work.  We will give a definition of the signs in terms of spectral
flow, as in Taubes' work~\cite{taubes:casson}.


The definition, slightly extended to the case of $SO(3)$--bundles, is
as follows. Fix a $w \in H^2(\Sigma';\Z_2)$ and consider the moduli
space $\R^*_w (\Sigma',SO(3))$ of the (gauge equivalence classes of)
irreducible flat connections in an $SO(3)$--bundle $P'_w$ over $\Sigma'$
having $w_2(P'_w) = w$. We say that $\R^*_w (\Sigma',SO(3))$ is
\emph{non-degenerate} if, for every $\alpha'\in \R^*_w (\Sigma',SO(3))$,
the cohomology $H^1(\Sigma';\ad\alpha')$ vanishes. If $\R^*_w(\Sigma',
SO(3))$ is non-degenerate, define the \emph{Boyer-Nicas type invariant}
\begin{equation}\label{E:bn}
\lambda_w (\Sigma') = \frac 1 2\; \sum_{\alpha'\in\R^*_w(\Sigma',SO(3))}
(-1)^{\mu(\alpha')},
\end{equation}
where $\mu (\alpha')$ is the spectral flow of the odd signature operator
(\ref{E:ka}) along a generic path connecting $\theta_w$ to $\alpha'$.
Here, $\theta_0$ is the product connection in the trivial bundle $P'_0$,
and $\theta_1$ is the (non-trivial) flat reducible connection in $P'_1$
with holonomy
\[
\begin{CD}
\pi_1 (\Sigma') @>>> H_1 (\Sigma';\Z_2) = \Z_2 @>\beta >> SO(3),
\end{CD}
\]

where $\beta$ maps the generator of $\Z_2$ to the diagonal matrix $\diag
(1,-1,-1) \in SO(3)$.

If $\R^*_w (\Sigma', SO(3))$ fails to be non-degenerate, it needs to be
perturbed first. In this case, the invariant is best described in terms
of the moduli space of projectively flat $U(2)$--connections\,: roughly
speaking, $SO(3)$ perturbations are $H^1 (\Sigma'; \Z_2)$--equivariant
$U(2)$ perturbations. We give a brief outline of the construction below
and refer to \cite{ruberman-saveliev:casson} for details.

For each $w \in H^2 (\Sigma';\Z_2)$ choose a $U(2)$--bundle $\tilde P'_w$
whose associated $SO(3)$--bundle is $P'_w$  (such a bundle exists because
$H^3 (\Sigma';\Z)$ is torsion free). Then $c_1 (\tilde P'_w) = w\pmod 2$,
and we make our choice so that $\tilde P'_w$ is a trivial bundle whenever
$w = 0 \pmod 2$. Every connection on $\tilde P'_w$ induces and is induced
by unique connections on $P'_w$ and on the line bundle $\det\tilde P'_w$.

Let us fix a connection $C$ on $\det\tilde P'_w$ and let $\A_w(\Sigma',
U(2))$ be the space of connections on $\tilde P'_w$ compatible with $C$.
If $w = 0 \pmod 2$ we choose $C$ to be the trivial connection; note that
$C$ plays no geometric role as different choices lead to equivalent
theories. Let $\B_w (\Sigma',U(2))$ be the quotient space of $\A_w
(\Sigma',U(2))$ by the action of gauge group consisting of unitary
automorphisms of $\tilde P'_w$ of determinant one. Contained in $\B_w
(\Sigma',U(2))$ is the moduli space $\R_w(\Sigma',U(2))$ of projectively
flat connections in $\tilde P'_w$. There is a natural action of the
group $H^1 (\Sigma';\Z_2)$ on both of the above spaces, with quotients
$\B_w (\Sigma',SO(3))$ and $\R_w (\Sigma',SO(3))$, respectively.

\begin{lemma}\label{L:action}
The action of $H^1(\Sigma';\Z_2)$ on $\R^*_w (\Sigma',U(2))$ is free,
and its orbit space is $\R^*_w(\Sigma',SO(3))$.
\end{lemma}

\begin{proof}
If $n$ is odd the group $H^1 (\Sigma';\Z_2)$ vanishes therefore we may
assume from the beginning that $n$ is even so that $H^1 (\Sigma';\Z_2)
= \Z_2$.

First let $w = 0 \pmod 2$ then $\R^*_w (\Sigma',U(2)) = \R^* (\Sigma',
SU(2))$.  The action of $\chi \in H^1 (\Sigma';\Z_2)$ viewed as a
homomorphism $\chi\co  \pi_1 \Sigma' \to \Z_2 = \{\,\pm 1\,\}$ is given
by $\alpha'\mapsto \chi\cdot\alpha'$. If the conjugacy class of
$\alpha'$ is fixed by $\chi\ne 1$ then $\chi\alpha' = u\alpha'u^{-1}$
for some $u\in SU(2)$ such that $u^2 = -1$. This means that $\alpha'$
is a binary dihedral representation and hence its lift $\pi^*\alpha'\co 
\pi_1\Sigma \to SU(2)$ is reducible. The fact that $\Sigma$ is an
integral homology sphere implies that $\pi^*\alpha'$ is trivial,
which means that $\alpha'$ cannot be irreducible.

Now suppose that $w = 1 \pmod 2$ then $\R^*_w(\Sigma',U(2))$ consists
of irreducible projective representations $\alpha'\co  \pi_1 \Sigma' \to
SU(2)$ with $w_2(\ad\alpha') = 1\pmod 2$. Arguing as above, we conclude
that the fixed points of a non-trivial $\chi \in H^1 (\Sigma';\Z_2)$
are the binary dihedral projective representations. If $\alpha'$ is
binary dihedral, the image of $\pi^*(\ad\alpha') = \ad(\pi^*\alpha')\co 
\pi_1 \Sigma \to SO(3)$ belongs to $SO(2)$, which makes it a trivial
representation. Therefore, $\ad\alpha'$ factors through
$\pi_1\Sigma/\pi_1\Sigma' = \Z_n$ and $\alpha'$ cannot be irreducible.
\end{proof}

The perturbation theory as described in Section \ref{S:perturb} extends
naturally to the case of projectively flat connections, see for instance
\cite[Section 5]{ruberman-saveliev:casson}. For a generic admissible
perturbation $h$ of type (\ref{ha}), compare with \cite{taubes:casson},
the perturbed projectively flat moduli space $\R^*_{w,h} (\Sigma',U(2))$
is non-degenerate. In particular, it consists of finitely many points,
and we define the \emph{Boyer-Nicas type invariant} by the formula
\begin{equation}\label{E:bn'}
\lambda_w(\Sigma') = \frac {1}{2\,|H_1(\Sigma';\Z_2)|}\;
\sum_{\alpha'\in\R^*_{w,h}(\Sigma',U(2))}(-1)^{\mu(\alpha')},
\end{equation}
where $\mu (\alpha')$ is the spectral flow of the perturbed odd signature
operator (\ref{E:ka}) along a generic path in $\A_w(\Sigma', U(2))$ which
connects a $U(2)$--lift of $\theta_w$ to a lift of $\alpha'$.

According to Lemma \ref{L:action}, definitions (\ref{E:bn}) and
(\ref{E:bn'}) agree if $\R^*_w (\Sigma',SO(3))$ is non-degenerate.
Moreover, according to \cite[Corollary 5.5]{ruberman-saveliev:casson},
a small generic admissible perturbation $h$ which makes $\R^*_{w,h}
(\Sigma',U(2))$ non-degene\-rate can be chosen to be equivariant with
respect to the action of $H^1(\Sigma';\Z_2)$. If $h$ is small enough,
this action still has no fixed points on $\R^*_{w,h} (\Sigma', U(2))$.
In addition, it preserves the spectral flow modulo 2, see \cite[pages
239-240]{braam-donaldson:knots}, so $\mu(\alpha')$ in (\ref{E:bn'})
can be replaced by the spectral flow along a path of $SO(3)$--connections
from $\theta_w$ to $\alpha'$. This allows to view (\ref{E:bn'}) as one
half times the signed count of the perturbed flat $SO(3)$--connections.

\begin{theorem}\label{thm2}
The invariants $\lambda_w (\Sigma')$ are well defined, and $\lambda_0
(\Sigma')$ equals $n$ times the Boyer-Nicas invariant of $\Sigma$.
Moreover, the sum of the $\lambda_w (\Sigma')$ over $w\in H^2 (\Sigma';
\Z_2)$ equals the equivariant Casson invariant $\lambda^{\tau}(\Sigma)$.
\end{theorem}

\begin{proof}
Since $\Sigma'$ is cyclically finite, all reducible flat connections in
the flat moduli space are isolated. This also remains true after a small
admissible perturbation. Therefore, we can use the standard cobordism
argument to show the well-definedness. Namely, for a generic path of
small admissible perturbations $h(t)$,\, $t\in [0,1]$, the parameterized
moduli space
\[
\bigcup_{t\in [0,1]}\; \R^*_{w,h(t)} (\Sigma',U(2))\times\{\,t\,\}
\]
gives a \emph{compact} oriented cobordism between the perturbed flat
moduli spaces $\R^*_{w,h(0)}(\Sigma',U(2))$ and $\R^*_{w,h(1)}
(\Sigma',U(2))$.
This allows us to conclude that the signed count of points in $\R^*_{w,h}
(\Sigma',U(2))$ is independent of $h$. A minor modification of Taubes'
argument \cite{taubes:casson} and the above discussion of equivariant
perturbations show that $\lambda_0 (\Sigma')$ equals $n$ times the
Boyer-Nicas invariant of $\Sigma$.

To prove the final statement, we use the identification
(Proposition~\ref{P:pullback}) between the moduli spaces $\R_w^{\tau}
(\Sigma,SU(2))$ and $\R^*_w (\Sigma',SO(3))$. If necessary, both spaces
can be perturbed by using admissible perturbations (\ref{ha}) on
$\Sigma'$ and their pull-backs to $\Sigma$, so that the perturbed moduli
spaces still match (the argument is completely similar to that in
Section 3.8 of \cite{collin-saveliev:casson}). Note that, since 
$\Sigma'$ is cyclically
finite, the perturbations may be chosen so that they do not affect the
reducible solutions.

Given a perturbed flat $SO(3)$--connection $\alpha'$ in $P'_w$, choose
a generic path of connections $\alpha'(t)$ connecting $\theta_w$ to
$\alpha'$. Lift it to a path $\pi^* \alpha'(t)$ of $SU(2)$--connections
over $\Sigma$ equivariant with respect to a certain constant lift $u\co  P
\to P$ and connecting $\theta$ to $\pi^*\alpha'$. Then $\mu^{\tau}
(\pi^*\alpha') = \mu(\alpha')$ due to the commutativity of the following
diagram
\[
\begin{CD}
(\Omega^0 \oplus \Omega^1)_u(\Sigma,\ad P) @> K^u_{\pi^*\alpha'(t)} >>
(\Omega^0 \oplus \Omega^1)_u(\Sigma,\ad P) \\
@A\cong A\pi^* A @A\cong A\pi^* A \\
(\Omega^0 \oplus \Omega^1)(\Sigma',\ad P'_w) @> K_{\alpha'(t)} >>
(\Omega^0 \oplus \Omega^1)(\Sigma,\ad P'_w).
\end{CD}
\]
\end{proof}

\begin{remark}\label{R:bn}
According to Theorem \ref{thm2}, no \emph{$H^1(\Sigma';\Z_2)$--equivariant}
perturbation is needed to compute $\lambda_w(\Sigma')$ using definition
(\ref{E:bn'}) -- any sufficiently small generic admissible perturbation $h$
which makes $\R^*_{w,h}(\Sigma',U(2))$ non-degenerate will do.
\end{remark}


\subsection{The formula}\label{S:formula}
Recall first that any homology lens space $\Sigma'$ such that $H_1
(\Sigma';\Z) = \Z_n$ can be obtained by $(n/q)$-surgery along a knot $k$
in an integral homology sphere $Y$, where $q$ is relatively prime to $n$
and $0 < q < n$, see for instance \cite{boyer-lines}.

\begin{proposition}\label{P:surgery}
Suppose that $\Sigma' = Y + (n/q)\,\cdot\, k$. Denote by $\lambda
(Y)$ the (regular) Casson invariant of $Y$, and by $\Delta (t)$ the
Alexander polynomial of $k\subset Y$ normalized so that $\Delta(1)=
1$ and $\Delta (t) = \Delta (t^{-1})$. Then, if $n$ is even,
\[
\lambda_w (\Sigma') = \frac n 2\;\lambda (Y) + \frac 1 8
\sum_{ \substack{m = 0 \\ m =\,w\mspace{-15mu}\pmod2}}^{n-1}
\sign^{m/n} k \; + \,\frac q 4\;\Delta''(1),\quad w = 0, 1,
\]
and, if $n$ is odd,
\[
\lambda_0 (\Sigma') = n\,\lambda (Y) + \frac 1 8\,\sum_{m=0}^{n-1}
\;\sign^{m/n} k \; + \, \frac q 2\;\Delta''(1).
\]
\end{proposition}

By $\sign^{\alpha} k$ we mean the Tristram--Levine equivariant knot
signature, see Section \ref{S6}.  The following result is immediate
from Theorem \ref{thm2} and Proposition \ref{P:surgery}.

\begin{corollary}\label{C:surgery}
For any $n\ge 1$, the equivariant Casson invariant can be expressed
by the formula
\begin{equation}\label{E:lambda'}
\lambda^{\tau} (\Sigma) = n\,\lambda (Y) + \frac 1 8\,
\sum_{m = 0}^{n-1}\; \sign^{m/n} k + \frac q 2\, \Delta''(1).
\end{equation}
\end{corollary}

\begin{remark}
Boyer and Lines \cite{boyer-lines} inquired whether the invariant of the
homology lens space $\Sigma' = Y + (n/q)\cdot k$ defined by the formula
\[
\overline{\lambda}(\Sigma') = \lambda (Y) + \frac 1 {8n}\;
\sum_{m = 0}^{n-1}\; \sign^{m/n} k + \frac q {2n}\, \Delta''(1)
\]
coincides with the Boyer-Nicas invariant $(1/n)\cdot\lambda_0(\Sigma')$,
see Question 2.24 of \cite{boyer-lines}. Proposition \ref{P:surgery}
implies that indeed $n\cdot\overline{\lambda}(\Sigma') = \lambda_0
(\Sigma')$ if $n$ is odd; however, $n\cdot\overline{\lambda}(\Sigma') =
\lambda_0(\Sigma') + \lambda_1(\Sigma')$ if $n$ is even.
\end{remark}

The rest of this section is devoted to the proof of Proposition
\ref{P:surgery}.


\subsection{Small scale perturbations}\label{S:small}
\begin{figure}[ht!]\small
\centering
\psfrag{phi}{$\phi$}
\psfrag{psi}{$\psi$}
\psfrag{pi}{$\pi$}
\psfrag{pi/2}{$\pi/2$}
\psfrag{0}{$0$}
\includegraphics[width=1.8in]{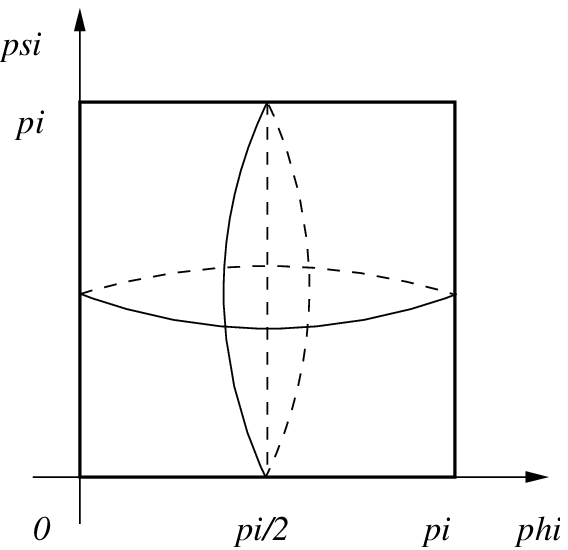}
\caption{}
\label{pillowcase}
\end{figure}

Let $X = Y \setminus N(k)$ be the knot $k$ exterior, and $m$ and $\ell$ be
the canonical meridian and longitude on $\p X$. Let $\P = \R(\p X,SU(2))$
be the $SU(2)$ pillowcase with coordinates $(\phi, \psi)$ such that the
holonomies along $m$ and $\ell$ equal $\exp(i\phi)$ and $\exp(i\psi)$,
respectively, shown in Figure \ref{pillowcase}. The inclusion $\p
X\to X$ induces a natural restriction map $r\co  \R(X,SU(2)) \to \P$ whose
image is generically one-dimensional.

\begin{figure}[ht!]\small
\centering
\psfrag{phi}{$\phi$}
\psfrag{psi}{$\psi$}
\psfrag{2pi}{$2\pi$}
\psfrag{0}{$0$}
\psfrag{C}{$\tilde\C$}
\includegraphics[width=1.8in]{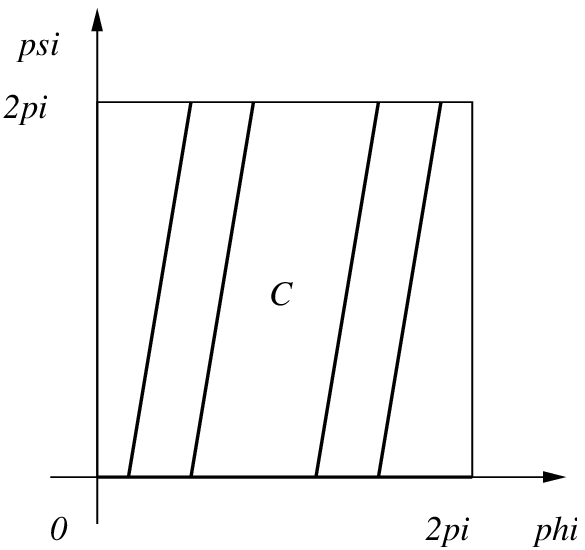}
\caption{}
\label{trefoilrep}
\end{figure}

Technically, it is more convenient to work with double covers of the above
as in \cite{herald:alexander}. More precisely, if we fix the maximal torus
$U(1)= \{\,\exp(i\phi)\,\}\subset SU(2)$ then the space of flat
$U(1)$--connections
on $\p X$ modulo $U(1)$--gauge group is a natural double cover $\tilde\P =
T^2$ of $\P$ branched at the central orbits. It gives rise to the double
cover $\tilde\R (X, SU(2))$ and the restriction map $r\co  \tilde\R (X,SU(2))
\to \tilde\P$. The 2-torus $\tilde\P$ will be depicted as the square
$\tilde\P = \{\,(\phi,\psi)\,|\,0\le\phi\le 2\pi,\, 0\le\psi\le 2\pi\,\}$
with properly identified sides; the double covering $\tilde\P\to \P$ then
corresponds to the quotient map by the action $\sigma (\phi,\psi) = (2\pi
- \phi, 2\pi - \psi)$ (which is geometrically a $180^{\circ}$ rotation of
the square $\tilde\P$ about the point $(\pi,\pi)$). The image of $r\co 
\tilde\R(X,SU(2))\to \tilde\P$ is invariant with respect to this action.
For instance, if $k$ is the left-handed trefoil, the image of $\tilde\R(X,
SU(2))$ in $\tilde\P$ consists of the circle $\{\,\psi = 0\,\}$ and two
open arcs limiting to it, as shown in bold in Figure \ref{trefoilrep}.

We will perturb the flat moduli space if necessary using admissible
perturbations $h$ described in \cite{herald:perturbations} (as amended
in \cite{herald:pa9}). The moduli space of $h$-perturbed flat
connections will be denoted by $\R_h (X,SU(2))$, and its double cover
by $\tilde\R_h(X,SU(2))$.

According to \cite{herald:perturbations}, for a generic admissible
perturbation $h$, the perturbed flat moduli space $\tilde\R_h(X,SU(2))$
consists of two central orbits, two smooth open arcs of abelian orbits
with one non-compact end limiting to each central orbit (together with
the central orbits, they form a circle), and a smooth 1-dimensional
manifold of irreducible orbits with finitely many ends, each limiting to
a different point on the abelian arc. If a point on the abelian arc with
coordinate $\phi$ is the limiting point of an irreducible arc then $\exp
(2i\phi)$ is a root of $\Delta(t)$. Since (see Section
\ref{S:rohlin}) the $n$--fold cyclic branched cover
of $k\subset Y$ is an integral homology sphere, $\Delta(\exp(2\pi i
m/n))\ne 0$ for $m \in \Z$, and we can conclude that

\begin{enumerate}
\item[(a)] no irreducible arc limits to a point on the abelian arc with
coordinate $\phi = \pi m/n$.
\end{enumerate}

Furthermore, the restriction map $r\co  \tilde\R_h (X, SU(2)) \to \tilde\P$,
which is well defined because $h$ is supported away from $N(k)$, is an
immersion taking the 1-dimensional strata of $\tilde\R_h (X,SU(2))$ into
the part of $\tilde\P$ away from the branching points. The map $r$
restricted to the circle of reducible orbits is an embedding with image
$\{\,\psi = 0\,\}$. The image of the irreducible part of $\tilde\R_h (X,
SU(2))$ in $\tilde\P$ will be denoted by $\tilde\C$ (compare with Figure
\ref{trefoilrep}). By using another small admissible perturbation if
necessary, one can ensure the following\,:

\begin{enumerate}
\item[(b)] the intersections of $\tilde\C$ with all of the circles
$\S (\pi m/n)\subset \tilde\P$ given by the equation $\phi = \pi
m/n$, $0\le m\le 2n$, are transverse (see \cite{herald:alexander}, Lemma 5.1);
\item[(c)] the intersection of $\tilde\C$ with the circle $\psi = \pi$
is transverse, and the $\phi$--coordinates of the intersection points
are different from $\pi m/n$, $0\le m\le 2n$. This can be achieved by
additional perturbations on loops $\gamma$ which are contained in a
collar neighborhood of $\p X$ and have the property that $\lk(\gamma,
k) \ne 0$ (see \cite{herald:perturbations}, page 408);
\item[(d)] the intersections of $\tilde\C$ with the curve   $\psi =
- (n/q)\,\cdot\,\phi$, and also with the curve $\psi = \pi - (n/q)\,
\cdot\,\phi$ if $n$ is even, are transversal (using the loops $\gamma$
as in (c) again). These intersections correspond to the respective
$\R^*_{w,h} (\Sigma',U(2))$, as we will explain below.
\end{enumerate}

\noindent
Note that, after all the admissible perturbations, $\tilde\C$ remains
invariant with respect to the involution $\sigma$.

\subsection{Large scale perturbations}
Recall that $\Sigma'$ is given by $(n/q)$ Dehn surgery along a knot
$k\subset Y$, ie, by adding a solid torus to $X$.  Parameterize this
solid torus by $\gamma\co  S^1\times D^2\to\Sigma'$, and write $k'$ for
its core $S^1 \times \{0\}$.  Let $f\co  SU(2)\to\mathbb R$
be a smooth function which is invariant with respect to conjugation.
Such a function is determined by a function $\tilde f\co  [0,\pi]\to\mathbb
R$ by the formula $\tilde f(t) = f(\exp(i t))$. Note that the smoothness
of $f$ requires that $\tilde f'$ vanish at $t = 0$ and $t = \pi$. We will
extend $\tilde f$ to a function on $[0,2\pi]$ by the rule $\tilde f(t) =
\tilde f(2 \pi - t)$ for $\pi\le t\le 2\pi$. Let $h_0$ be the admissible
perturbation given by the formula (compare with (\ref{ha}))
\begin{equation}\label{h0}
h_0(\alpha) = \int_{D^2}\; f(\hol_{\alpha}(\gamma (S^1\times\{z\}),s))
\,\eta(z)\,d^2 z.
\end{equation}
Notice that the knot exterior $X$ is naturally a subset of $\Sigma'$.
Any perturbed flat connection over $\p X$ is actually flat hence it
is determined by two parameters $a$, $b$ where, in a suitable
parametrization, $\exp(i a)$ is the holonomy around meridian (the
boundary of the disc $D^2$) and $\exp(i b)$ is the holonomy around
the longitude (a curve parallel to $k$). The following is a minor
modification of Lemma 4 of \cite{braam-donaldson:knots}.

\begin{lemma}
For sufficiently small $\ep > 0$, restriction to $X$ defines a one-to-one
correspondence between

\begin{itemize}\label{L:extend}
\item[\rm(a)] the gauge equivalence classes of $(h_0 + \ep h)$--perturbed
flat connections on $\Sigma'$, and
\item[\rm(b)] the gauge equivalence classes of $\ep h$--perturbed flat
connections on $X$ satisfying the equation $a = \tilde f'(b)$.
\end{itemize}
\end{lemma}

\begin{figure}[ht!]\small
\centering
\psfrag{f}{$\phi$}
\psfrag{s}{$\psi$}
\psfrag{2pi}{$2\pi$}
\psfrag{pi}{$\pi$}
\psfrag{1}{$\psi = -(n/q)\cdot\phi$}
\psfrag{b}{$a = \tilde f'(b)$}
\psfrag{0}{$0$}
\includegraphics{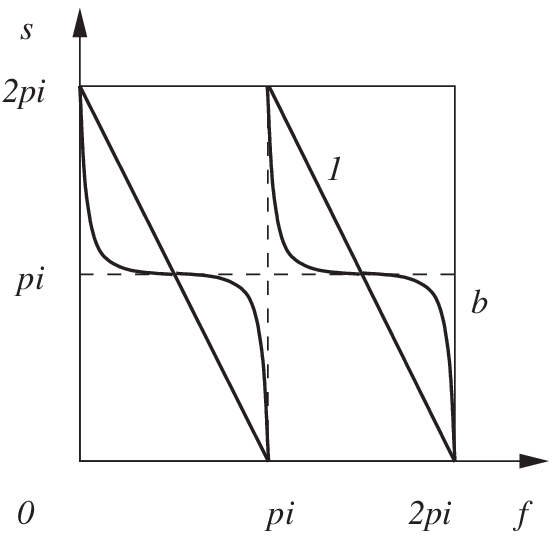}

\caption{\,:\, $n = 2$, $q = 1$}
\label{flarge}
\end{figure}

The extension condition (b) in Lemma~\ref{L:extend} (ie,  the
equation $a = \tilde f'(b)$) defines a curve in $\tilde \P$.  In the
case that $f\equiv 0$, this curve would be the line $\psi = -
(n/q)\,\cdot\,\phi$, using coordinates $(\phi,\psi)$ dual to the
meridian and longitude of $k$ (which differ from those of $k'$).  A
more realistic example may look as shown in Figure \ref{flarge}.
The figure shows, for $(n/q) = (2/1)$, the curves
$\psi = - (n/q)\,\cdot\,\phi$ and $a = \tilde f' (b)$.  Note that the
curve $a = \tilde f'(b)$ is invariant with respect to $\sigma$.

\subsection{The case of even $n$}
Let us place $\tilde\C$ in general position as in Section \ref{S:small},
then $\R^*_{0,h}(\Sigma',U(2))$ is in one-to-one correspondence with the
orbits of $\sigma$ acting on $\tilde\C\,\cap\,\{\,\psi= -(n/q)\cdot\phi
\,\}$. Note that this action is free, see Lemma \ref{L:action}, so each
such orbit has two elements.  Taking into account the $1/4$
factor in the definition~\eqref{E:bn'} and the orientations set in
Section \ref{S:orient}, the invariant $\lambda_0(\Sigma')$ equals one
eighth times the intersection number of $\tilde\C$ with the curve
$\psi = -(n/q)\,\cdot\,\phi$. This can be proved along the same lines as
\cite[Theorem 7.1]{herald:alexander}. The heart of that proof is a
spectral flow calculation that does not require that $\Sigma'$ be an
integral homology sphere.

The curve $\psi= -(n/q)\,\cdot\,\phi$ wraps $n$ times around $\tilde \P$
in the $\psi$--direction, and the points of intersection with the curve
$\psi = 0$ separate it into $n$ segments given by $2 \pi jq/n \le \phi
\le 2\pi (j+1)q/n$ where $j = 0,\ldots, n-1$. By choosing a suitable
function $f$, one can deform each of these segments in a
$\sigma$--equivariant manner into a graph $a = \tilde f'(b)$ which
lies arbitrarily close to the union
\[
\tilde S (2\pi jq/n)\,\cup\,\tilde S (2\pi(j+1)q/n)\,\cup\,
\{\,\psi = \pi\,\},
\]
as shown in Figure \ref{flarge}. With orientations set as in Section
\ref{S:orient}, we see that the intersection number of $\tilde\C$ with
the curve $\psi = -(n/q)\,\cdot\,\phi$ is the sum of two terms. The
first one is the intersection number of $\tilde\C$ with the union of
the circles $\tilde S (2\pi j/n)$ over $j = 0,\ldots, n-1$. According
to \cite{herald:alexander}, it equals
\[
4n\cdot\lambda(Y) + \frac 1 2\,\sum_{j=0}^{n-1}\;\sign^{2j/n} k = 4n
\cdot\lambda(Y) + \sum_{\substack {m=0 \\ m\,\text{even}}}^{n-1}\;
\sign^{m/n} k.
\]
The second term equals $q$ times the intersection number of  $\tilde\C$
with the circle $\{\,\psi=\pi\,\}$; this is given by $-2\,\Delta''(1)$.
Up to an overall sign, this can be seen as in \cite[Proposition 7.2]
{herald:alexander} using the gauge theoretic description of $\Delta''(1)$
given in \cite{masataka}; the sign is then fixed with the help of the Casson
surgery formula $\lambda(Y - k) = \lambda(Y) - (1/2)\,\Delta''(1)$.
This proves the formula for $\lambda_0(\Sigma')$.

\begin{figure}[ht!]\small
\centering
\psfrag{p}{$\pi/2$}
\psfrag{1}{$\theta_1$}
\psfrag{phi}{$\phi$}
\psfrag{s}{$\psi$}
\psfrag{2}{$2\pi$}
\psfrag{0}{$0$}
\psfrag{q}{$\psi =\pi -(n/q)\cdot\phi$}
\includegraphics{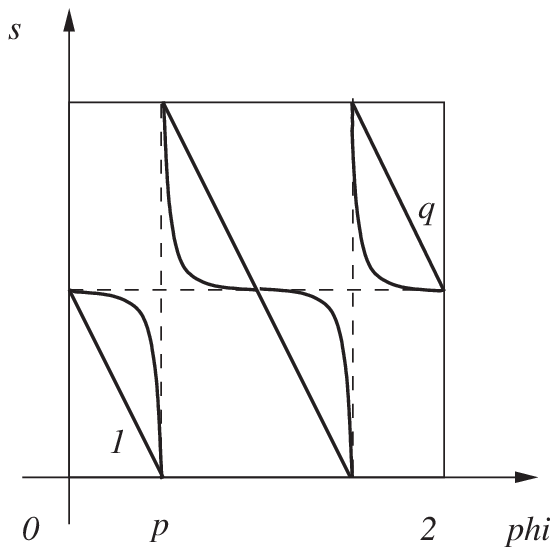}
\caption{\,:\, $n = 2$, $q = 1$}
\label{theta1}
\end{figure}

To prove the formula for $\lambda_1 (\Sigma')$, we consider the
intersection $\tilde\C\,\cap\,\{\,\psi = \pi - (n/q)\,\cdot\,\phi\,\}$
and deform the curve $\psi = \pi - (n/q)\,\cdot\,\phi$ again as shown
in Figure \ref{theta1}. The rest of the proof follows the lines of the
above proof for $\lambda_0 (\Sigma')$. To express $\lambda_1(\Sigma')$
as the intersection number in $\tilde \P$, we again refer to the proof
of \cite[Theorem 7.1]{herald:alexander}, where the trivial connection
should be replaced with $\theta_1$, see Figure \ref{theta1}.

\subsection{The case of odd $n$}
The proof is completely similar to that for $\lambda_0 (\Sigma')$ in the
even case. The invariant $\lambda_0 (\Sigma')$ equals one fourth times the
intersection number of $\tilde\C$ with the curve $\psi = -(n/q)\,\cdot\,
\phi$ (after fixing orientations). Deforming $\psi = -(n/q)\,\cdot\,\phi$
as shown in Figure \ref{flarge}, we see that the above intersection number
equals the sum of two terms. The first one is the intersection number of
$\tilde\C$ with the union of the circles $\tilde S (2\pi j/n)$ over $j =
0,\ldots, n-1$. According to \cite{herald:alexander}, it equals
\[
4n\cdot\lambda(Y) + \frac 1 2\,\sum_{j=0}^{n-1}\;\sign^{2j/n} k = 4n
\cdot\lambda(Y) + \frac 1 2\,\sum_{m = 0}^{n-1}\;\sign^{m/n} k
\]
(we used the fact that $\sign^{\alpha} k = \sign^{1 + \alpha} k$).
The second term as before equals $q$ times the intersection number
of $\tilde\C$ with the circle $\{\,\psi =\pi\,\}$, which equals $-2\,
\Delta''(1)$. This completes the proof.


\subsection{Orientations}\label{S:orient}
The tangent space of $\tilde \P$ at any point is canonically isomorphic
to $H^1(\p X; \mathbb R)$ and hence is canonically oriented after we
orient $\p X$ as the boundary of $X$. We adapt the convention that an
outward normal vector followed by an oriented basis for $T_x (\p X)$
gives an oriented basis for $T_x X$ (this is opposite to the convention
used in \cite{herald:alexander}).

Orientations of various curves in the above construction depend on
an orientation of $H^1 (X; \mathbb R)$ or, equivalently, a choice of
meridian. Once it is fixed, the circles $\{\,\psi = 0\,\}$ and $\{\,
\psi = \pi\,\}$ are naturally oriented by the parameter $\phi$. The
circles $S(\alpha)$ are oriented by the requirement that the
intersection number of $\{\,\psi = 0\,\}$ with $S(\alpha)$ be one;
the same requirement orients the curves $\{\,\psi = -(n/q)\cdot\phi
\,\}$ and $\{\,\psi = \pi - (n/q)\cdot\phi\,\}$. The curve $\tilde\C$
is oriented as described in \cite{herald:alexander}, Section 2.

Note that $\sigma$ reverses the orientation of $H^1 (X;\mathbb R)$;
however, it preserves all the intersection numbers.


\subsection{Examples}
Let $p$, $q$, and $r$ be pairwise relatively prime positive integers,
and consider the Brieskorn homology sphere
\[
\Sigma(p,q,r) = \{\,x^p+y^q+z^r=0\,\}\,\cap\,S^5\subset\mathbb C^3.
\]
The map $\tau(x,y,z) = (x,y,e^{2\pi i/r}z)$ defines a $\Z_r$--action
on $\Sigma(p,q,r)$ whose fixed point set is the singular fiber  $k_r
= \Sigma(p,q,r)\,\cap\,\{\,z=0\,\}$. The quotient of $\Sigma(p,q,r)$
by this action is $S^3$, with branch set the right--handed
$(p,q)$--torus knot $k$. In particular, if we write $\Sigma(p,q,r) =
X\cup N(k_r)$ where $X$ is the knot exterior, then the action on
$N(k_r) = S^1\times D^2$ will be given by the formula $(s,z) \mapsto
(s,e^{2\pi i/r}z)$.

The Brieskorn homology sphere $\Sigma(p,q,r+pq)$, obtained from
$\Sigma(p,q,r)$ by $(-1)$--surgery along $k_r$, admits a \emph{free}
$\Z_r$--action which extends the action on $X$ to $S^1\times D^2$
by the formula $(s,z)\mapsto (e^{2\pi i/r}s, e^{2\pi i/r}z)$. If
$\Sigma(p,q,r+pq)$ is viewed as a link of singularity, this action
is given by
\begin{equation}\label{E:action}
\tau(x,y,z) = (e^{2\pi iq/r}x,e^{2\pi ip/r}y,e^{2\pi i/r}z).
\end{equation}
Corollary \ref{C:surgery} can be used to compute the equivariant
Casson invariant of $\Sigma(p,q,r+pq)$ with respect to the action
(\ref{E:action}) as follows.

Let $\Sigma'$ be the quotient manifold of $\Sigma(p,q,r+pq)$. It
is immediate from the above description that $\Sigma'$ can be
obtained from $S^3$ by $(-r)$--surgery along the $(p,q)$--torus
knot $k$. Therefore,
\begin{equation}\label{E:seifert}
\lambda^{\tau}(\Sigma(p,q,r+pq)) = \frac 1 8 \sum_{m=0}^{r-1}
\sign^{m/r} k - \frac 1 2\,\Delta''_k (1),
\end{equation}
where
\begin{alignat}{2}
\frac 1 8 \sum_{m=0}^{r-1}\sign^{m/r} k
&= \lambda^{\tau}(\Sigma(p,q,r)), & &\quad\text{by (\ref{E:cs}),}
\notag \\
&= \lambda (\Sigma(p,q,r)), & &\quad\text{according to
\cite{collin-saveliev:casson},}
\notag
\end{alignat}

and $\Delta''_k (1) = \Delta''_{k_r} (1) = (p^2 - 1)(q^2 - 1)/12$,
see \cite{neumann-wahl}. Now, the right hand side of (\ref{E:seifert})
equals $\lambda(\Sigma(p,q,r)) - (1/2)\,\Delta''_{k_r}(1)$, which
coincides with $\lambda(\Sigma(p,q,r+pq))$ by Casson's surgery
formula. Therefore,
\[
\lambda^{\tau}(\Sigma(p,q,r+pq)) = \lambda(\Sigma(p,q,r+pq)).
\]
Another way to see this would be using the fact that $\tau$ can be
included into the natural $S^1$--action on the link $\Sigma(p,q,r+
pq)$ given by
\[
(x,y,z) \mapsto (t^{q(r+pq)}x, t^{p(r+pq)}y, t^{pq}z)
\]
and is therefore isotopic to the identity.

\section{The Rohlin invariant}\label{S:rohlin}
Our next goal is prove that $\lambda^{\tau}(\Sigma)$ equals the Rohlin
invariant $\rho(\Sigma)$ modulo 2, and thus verify Theorem~\ref{T:reduction}.
This is a known fact when $\Fix(\tau)\ne\emptyset$,
see~\cite{collin-saveliev:casson}. In the case of $\Fix(\tau) = \emptyset$,
we will achieve our goal by identifying the right hand side of
(\ref{E:lambda'}) with $\rho(\Sigma)$ modulo 2.

\subsection{Reduction to the Arf--invariants}
We use notations of Section \ref{S:formula}.
Let $Y_n$ be the $n$--fold branched cover of $Y$ with branch set the knot
$k$; the preimage of $k$ in $Y_n$ is a knot which we call $k_n$.  Observe
that $\Sigma = Y_n + (1/q)\,\cdot k_n$ so that $Y_n$ is an integral
homology sphere. According to (\ref{E:cs}),
\[
n\cdot\lambda (Y) + \frac 1 8\,\sum_{m = 0}^{n-1}\; \sign^{m/n} k =
\lambda^{\tau} (Y_n)
\]
is the equivariant Casson invariant of $Y_n$, whose modulo 2 reduction
is known~\cite[Theorem 2]{collin-saveliev:casson} to equal $\rho(Y_n)$.
Since $(1/2)\,
\Delta''(1) = \arf(k)\pmod 2$, we only need to show that $\rho (\Sigma)
= \rho (Y_n) + q\cdot \arf (k)\pmod 2$.   Now, $\rho(\Sigma) = \rho(Y_n)
+ q\cdot\arf(k_n)\pmod 2$ by Rohlin's surgery formula. Therefore, the
proof will be complete after we show the following.

\begin{proposition}\label{P:arf}
Let $Y$ be an integral homology sphere and $\pi\co  Y_n \to Y$ be its
$n$--fold cyclic branched cover with branch set a knot $k$. Let
$k_n$ be the knot $\pi^{-1} (k)$ in $Y_n$. If $Y_n$ is an integral
homology sphere then $\arf (k_n) = \arf (k) \pmod 2$.
\end{proposition}

\subsection{An algebraic lemma}
Proposition \ref{P:arf} will be deduced from a purely algebraic lemma, using
Levine's theorem~\cite{levine:polynomials} connecting the Arf invariant
and the Alexander polynomial.  A polynomial $\Delta(t)\in \Z[t,t^{-1}]$ is
said to be \emph{symmetric and normalized} if $\Delta(t) = \Delta(t^{-1})$
and $\Delta(1) = 1$.  Note that this implies that $\Delta(-1) \equiv 1$ or
$5 \pmod 8$.  Levine showed that if the Alexander polynomial of a knot is
symmetric and normalized, the Arf invariant of the knot is $0$ if
$\Delta(-1)\equiv 1\pmod{8}$ and is $1$ if $\Delta(-1) \equiv 5 \pmod 8$.

\begin{lemma}\label{L:galois} Suppose that $\Delta(t), D(t) \in
\Z[t,t^{-1}]$ are symmetric, normalized polynomials.   For an integer
$n$, denote by $\omega$ a primitive $n^{th}$ root of unity.  Assume
that $D(t^n) = \prod_{j=0}^{n-1}\Delta(\omega^j t)$.  If $n$ is even,
then $D(-1) = \Delta(-1)\equiv 1 \pmod{8}$.  If $n$ is odd, then
$D(-1) = a\cdot\Delta(-1)$ where $a \in \Z$ is congruent to $1\pmod8$.
\end{lemma}

\begin{proof}  Write $n=2p+1$ if $n$ is odd, and $n=2p$ if $n$ is even.
If $\xi$ is a primitive $2n^{th}$ root of unity, then $\xi^n = -1$.
Write $\nu=\Delta(\xi)\Delta(\xi^3)\cdots\Delta(\xi^{2p-1})$.  Note
that  the symmetry of $\Delta$ means that $\Delta(\xi^k) =
\Delta(\xi^{2n-k})$.  If $n=2p$, this means that the last term in
$\nu$ is $\Delta(\xi^{4p-2p+1}) = \Delta(\xi^{2p+1})$.  Likewise, if
$n = 2p+1$, then the last term is $\Delta(\xi^{2p+3})$. Now compute
\begin{align*}
D(-1) &=\prod_{j=0}^{n-1}\Delta(\xi^{2j+1}) =
\Delta(\xi)\Delta(\xi^3)\cdots\Delta(\xi^{2p-1})\Delta(\xi^{2p+1})
\cdots\Delta(\xi^{2n-1}) \\
&=
\begin{cases}
\; \nu^2, & n\ {\rm even,}\\
\; \nu^2 \Delta(\xi^{2p+1}), &n\ {\rm odd.}
\end{cases}
\end{align*}
Now $D(-1)$, being congruent to $D(1) \pmod2$ is odd, as is
$\Delta(-1)$. We claim that $\nu$, which is \emph{a priori} an
algebraic integer, is actually an integer.  In the case that $n=2p+1$
is odd, $\Delta(\xi^{2p+1})= \Delta(-1)$, so $\nu^2$ will necessarily
be the square of an odd integer, and therefore congruent to $1
\pmod8$.  In the case that $n$ is even, we get by the same argument
that $D(-1)$ is congruent to $1\pmod8$.

To see that $\nu$ is an integer, we use a little Galois theory.  By
definition, $\nu \in \Z(\xi)$, which is the ring of integers in the
cyclotomic field $\Q(\xi)$. We need to show that $\nu$ is fixed by
every element of the Galois group of the field extension $\Q(\xi)/\Q$.
This group is the (multiplicative) group of units $\Z_{2n}^\times$,
where $r \in \Z_{2n}^\times$ acts by $\xi^j \to \xi^{jr}$.
Note that the set of  odd numbers $\{1,3,\dots,2n-1\} \subset \Z_{2n}$
is invariant under $j \to -j \pmod{2n}$.   If $n$ is even, the
orbits all have two elements, one $\le 2p-1$ and one $\ge 2p+1$. If
$n$ is odd, there is one orbit, namely $2p+1$, with one element while
all the others have two elements. Consider
the action of $r \in \Z_{2n}^\times$ on the exponents appearing in
$\nu$, namely $E:=\{1,3,\dots,2p-1\}$; note that $r$ preserves orbits
since $r(-a) = -ra$.  Now if $ra\in E $ then $-ra \not\in E$, and so
$rE$ consists of one element of each orbit (except $2p+1$ if $n$ is
odd).  It follows that $r\cdot \nu = r\cdot\prod_{a\in
E}\Delta(\xi^a) =\prod_{a\in E}\Delta(\xi^{ra}) =\prod_{a\in
E}\Delta(\xi^{\pm a}) = \nu$.

For odd $n$, we are finished, but for $n$ even, we still have to show
that $\Delta(-1) \equiv 1 \pmod8$.  In fact, $\Delta(-1) =1$, by the
following argument\,:
\begin{align*}\notag
D(1) &= \prod_{j=0}^{n-1}\Delta(\omega^{j+1}) \\
&= \Delta(1)\Delta(\omega^p) \Delta(\omega)\Delta(\omega^2) \cdots
\Delta(\omega^{p-1})\Delta(\omega^{p+1})\cdots\Delta(\omega^{n-1}) \\
&= \Delta(1)\Delta(-1)\eta
\end{align*}
for some $\eta$. Now the Galois group of $\Q(\omega)/\Q$ permutes the
elements $\omega^j$, and no element $\omega^j$ for $j \ne p$ gets sent
to $-1$. So the Galois group preserves $\eta$. Hence $\eta \in \Z$, so
$\Delta(-1) | D(1)$.  Since $D(1) =1$ and $\Delta$ is symmetric and
normalized, it follows that $\Delta(-1) =1$.
\end{proof}

\begin{remark}
The final argument is clearly equivalent to (and is inspired by) the fact
that there is a branched cover $Y_{2p} \to Y_2$, so the hypothesis
that $Y_n$ be a homology sphere implies that $Y_2$ is one also.
\end{remark}

\subsection{Proof of Proposition~\ref{P:arf}}
If $k$ is a knot in the homology sphere $Y$, and the $n$-fold cyclic
branched cover $Y_n$ is a homology sphere, then
the Alexander polynomials of $k$ and $k_n$ are related as in
Lemma~\ref{L:galois}, see \cite{fox:periodic}. Levine's theorem implies that
their Arf invariants must be the same.
\begin{remark}
Chuck Livingston has pointed out a more geometric proof of
Proposition~\ref{P:arf} that has the additional advantage of applying
to $\Z_2$ homology spheres.  Here is an outline of the argument.
Livingston observes that the Arf invariant of a knot in a $\Z_2$
homology sphere $Y$ can be defined as the spin-bordism class of a
Seifert surface for the knot, using the unique spin structure induced
from $Y$.  This spin structure, restricted to the knot complement,
pulls back under the branched covering, and clearly gives the same
spin structure on a lift of the Seifert surface.  By obstruction
theory, if the branched cover is also a a $\Z_2$ homology sphere,
this lifted spin structure (or one that is identical in a
neighborhood of the lifted Seifert surface) can be extended over the
branched covering.  Hence the Arf invariant upstairs is the same as
the Arf invariant downstairs.
\end{remark}


\section{Applications}\label{S:examples}
The relations between the equivariant Casson invariant on one hand, and
Donaldson polynomials and Floer Lefschetz numbers, on the other, lead
to several interesting applications which we describe in this section.

\subsection{Invariant $\lambda^{\tau}(\Sigma)$ as a function of $\tau$}
Let $\Sigma$ be an oriented integral homology sphere, and $\tau\co  \Sigma
\to \Sigma$ be an orientation preserving diffeomorphism of finite order.
Theorem \ref{T:finite} identifies $\lambda^{\tau}(\Sigma)$ with one
fourth of the degree--zero Donaldson polynomial of the mapping torus of
$\tau$.

If $\tau$ is (smoothly) isotopic to the identity, then the mapping torus
of $\tau$ is diffeomorphic to the product $S^1\times\Sigma$. Since the
degree--zero Donaldson polynomial is an invariant of diffeomorphism, we
conclude that $\lambda^{\tau} (\Sigma) = \lambda(\Sigma)$. It is not
clear if this result can be obtained directly from the definition of
$\lambda^{\tau} (\Sigma)$.  In the case that $\tau$ is contained in a
circle action, this was shown by direct calculation
in~\cite{collin-saveliev:casson}.

\subsection{Plumbed manifolds}
The equivariant Casson approach to the Floer Lefschetz number allows one to
actually compute that number in several interesting examples.

Given a plumbed homology sphere $\Sigma$, there is a preferred class of
orientation preserving involutions $\tau$ which represent $\Sigma$ as
a double branched cover of $S^3$ with branch set a Montesinos knot
$k_{\tau}$, see \cite{siebenmann:rohlin}. The knots corresponding to
different involutions from this class are related by a sequence of mutations.
Since the knot signature is a mutation invariant, this construction
ends up in a well defined integer valued invariant $\bar\mu(\Sigma) =
(1/8)\cdot\sign(k_{\tau})$ called the \emph{$\bar\mu$--invariant}, see
\cite{neumann:plumbing} and \cite{siebenmann:rohlin}. It is immediate from
(\ref{E:cs}) that $\bar\mu(\Sigma) = \lambda^{\tau}(\Sigma)$.

More can be said if $\Sigma$ is a Seifert fibered homology sphere
$\Sigma(a_1,\ldots,a_n)$. In this case, $\tau$ is induced by complex
conjugation on $\Sigma(a_1,\ldots,a_n)$ viewed as a link of complex
singularity. According to \cite{saveliev:seifert}, there exists an
equivariant admissible perturbation $h$ such that assumption $(*)$ of
Section \ref{S:lef} holds (whether the same holds for arbitrary
plumbed homology spheres remains to be seen). Therefore, we have the
following result.

\begin{proposition}
The $\bar\mu$--invariant of a Seifert homology sphere
$\Sigma(a_1,\ldots,a_n)$ equals half the Lefschetz number of $\tau_*\co 
I_*(\Sigma(a_1,\ldots,a_n)) \to I_*(\Sigma(a_1,\ldots,a_n))$.
\end{proposition}

According to \cite{fs:instanton}, the group $I_k(\Sigma(a_1,\ldots,a_n))$
is trivial if $k$ is even, and is free abelian group of rank, say, $b_k$
if $k$ is odd.

\begin{proposition}
The map $\tau_*\co I_k (\Sigma(a_1,\ldots,a_n)) \to I_k (\Sigma(a_1,
\ldots,a_n))$ is identity if $k = 1\pmod 4$ and minus identity if
$k = 3\pmod 4$. In particular, the Floer Lefschetz number of $\tau_*$
equals $-b_1 + b_3 - b_5 + b_7$.
\end{proposition}

\proof
The generators $\alpha$ of the Floer chain complex of $\Sigma (a_1,
\ldots, a_n)$ can be chosen so that all of them have odd Floer index
and are preserved by the action induced by $\tau$, see
\cite{saveliev:seifert}. According to \cite{saveliev:brieskorn} and
\cite{saveliev:seifert}, the equivariant spectral
flow for these generators is given by the formula $\mu^{\tau}(\alpha)
= \frac 1 2\,(\mu(\alpha) + 1)\pmod 4$. The result is now immediate
from the following calculation, compare with Section \ref{S:lef}\,:
$$
\tau_*(\alpha) = (-1)^{\mu^{\tau}(\alpha)-\mu(\alpha)}\,\alpha
= (-1)^{(k-1)/2}\,\alpha.
\eqno{\qed}$$

\subsection{Akbulut cork}
\begin{figure}[ht!]\small
\centering
\psfrag{0}{$0$}
\includegraphics{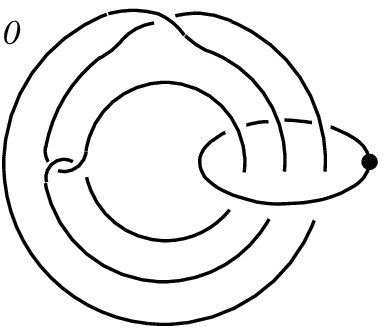}
\caption{}
\label{cork}
\end{figure}

Another situation when the Floer Lefschetz number can be computed
explicitly is that of the Akbulut cork. By \emph{Akbulut cork} we
mean the smooth contractible 4--manifold $W$ obtained by attaching
a two-handle to $S^1\times D^3$ along its boundary as shown in
Figure \ref{cork}. It can be embedded into a blown up elliptic surface
$E(n)\#(-\mathbb C P^2)$ in such a way that cutting it out and re-gluing
by an involution on $\Sigma = \p W$ changes the smooth structure on $E(n)
\# (-\mathbb C P^2)$ but preserves its homeomorphism type, see
\cite{akbulut:contractible} and \cite{gompf-stipsicz:book}.

\begin{figure}[ht!]\small
     \begin{minipage}[b]{0.50\linewidth}
      \centering
      \psfrag{0}{$0$}
      \includegraphics{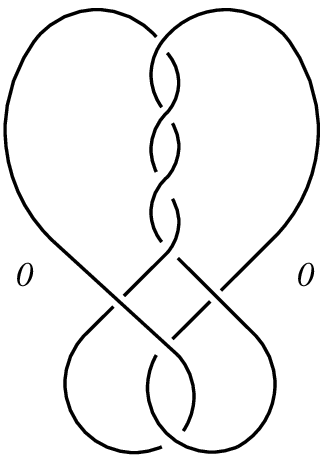}
\caption{}
\label{symmcork}
     \end{minipage}
     \begin{minipage}[b]{0.49\linewidth}
      \centering
      \psfrag{k}{$k^*$}
      \includegraphics{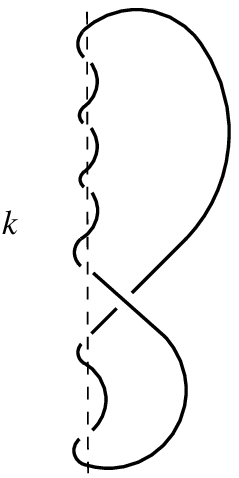}
\caption{}
\label{corkquotient}
     \end{minipage}
\end{figure}

The involution $\tau\co\Sigma \to \Sigma$ simply interchanges the two link
components, $k_1$ and $k_2$; this is best seen when the link is drawn in
a symmetric form as in Figure \ref{symmcork}. Note that $\Sigma$ is a
hyperbolic homology sphere.

The Floer homology of $\Sigma$ is trivial in even degrees, and is a copy
of $\Z$ in each of the odd degrees, see \cite{saveliev:cork}. Therefore,
$\tau_*\co  I_n (\Sigma) \to I_n (\Sigma)$ is necessarily plus or minus
identity for each $n$. It turns out that the exact answer can be obtained
from computing the Lefschetz number of $\tau_*$ keeping in mind that
condition $(*)$ of Section \ref{S:lef} is satisfied because the character
variety $\R^*(\Sigma)$ is non-degenerate, see \cite{saveliev:cork}.

The manifold $\Sigma/\tau$ is obtained from $S^3$ by surgery on the knot
$k^*$ which is the image of the link $k_1 \cup k_2$, see Figure
\ref{corkquotient}. Notice that the canonical longitudes of $k_1$ and $k_2$
project onto a longitude of $k^*$ whose linking number with $k^*$ equals
one. The dotted line in Figure \ref{corkquotient} represents the branch set
$k$.

This picture can be viewed as a surgery description of the knot $k$ in
$\Sigma/\tau = S^3$. A little exercise in Kirby calculus shows that $k$
is obtained from the left-handed $(5,6)$--torus knot on six strings by
adding one full left-handed twist on two adjacent strings. The signature
of $k$ can differ by at most two from the signature of the left-handed
$(5,6)$-torus knot, which equals 16. Since $\sign k$ must be divisible
by eight, we conclude that $\sign k = 16$. Therefore the Lefschetz number
of $\tau_*$ equals $2$. This implies that $\tau_* = -\id\co  I_* (\Sigma) \to
I_* (\Sigma)$.


\begin{thebibliography}

\bibitem{akbulut:contractible}
\textbf{Selman Akbulut}, \emph{A fake compact contractible {$4$}-manifold},
J. Differential Geom. 33 (1991) 335--356

\bibitem{boden-nicas}
\textbf{Hans~U Boden}, \textbf{Andrew Nicas}, \emph{Universal formulae for
{${\rm SU}(n)$} {C}asson invariants of knots}, Trans. Amer. Math. Soc. 352
(2000) 3149--3187

\bibitem{boyer-nicas}
\textbf{Steven Boyer}, \textbf{Andrew Nicas}, \emph{Varieties of group
representations and {C}asson's invariant for rational homology {$3$}-spheres},
Trans. Amer. Math. Soc. 322 (1990) 507--522

\bibitem{boyer-lines}
\textbf{Steven Boyer}, \textbf{Daniel Lines}, \emph{Surgery formulae for
{C}asson's invariant and extensions to homology lens spaces}, J. Reine
Angew. Math. 405 (1990) 181--220

\bibitem{braam-donaldson:knots}
\textbf{Peter~J Braam}, \textbf{Simon~K Donaldson}, \emph{Floer's work on
instanton homology, knots and surgery}, from: ``The Floer memorial volume'',
Progr. Math. 133, Birkh\"auser, Basel (1995)  195--256

\bibitem{braam-matic:smith}
\textbf{Peter~J Braam}, \textbf{Gordana Mati{\'c}}, \emph{The {S}mith
conjecture in dimension four and equivariant gauge theory}, Forum Math. 5
(1993) 299--311

\bibitem{collin-saveliev:casson}
\textbf{Olivier Collin}, \textbf{Nikolai Saveliev}, \emph{Equivariant
{C}asson invariants via gauge theory}, J. Reine Angew. Math. 541 (2001)
143--169

\bibitem{donaldson:orientation}
\textbf{Simon~K Donaldson}, \emph{The orientation of {Yang--Mills} moduli
spaces and 4-manifold topology}, J. Differential Geom. 26 (1987) 397--428

\bibitem{donaldson:floer}
\textbf{Simon~K Donaldson}, \emph{{Floer Homology Groups in Yang--Mills
Theory}}, Cambridge University Press (2002)

\bibitem{fs:instanton}
\textbf{Ronald Fintushel}, \textbf{Ronald J Stern}, \emph{Instanton
homology of Seifert fibred homology three spheres}, Proc. London Math.
Soc. 61 (1990) 109--137

\bibitem{fox:periodic}
\textbf{Ralph~H Fox}, \emph{Knots and periodic transformations}, from:
``Topology of 3-manifolds and related topics (Proc. The Univ. of Georgia
Institute, 1961)'', Prentice-Hall, Englewood Cliffs, NJ (1962)  177--182

\bibitem{frohman-nicas}
\textbf{Charles Frohman}, \textbf{Andrew Nicas}, \emph{An intersection homology
invariant for knots in a rational homology {$3$}-sphere}, Topology 33 (1994)
123--158

\bibitem{furuta-ohta}
\textbf{Mikio Furuta}, \textbf{Hiroshi Ohta}, \emph{Differentiable structures
on punctured $4$-manifolds}, Topology Appl. 51 (1993) 291--301

\bibitem{gompf-stipsicz:book}
\textbf{Robert~E Gompf}, \textbf{Andr{\'a}s~I Stipsicz}, \emph{$4$-manifolds
and {K}irby calculus}, American Mathematical Society, Providence, RI (1999)

\bibitem{herald:alexander}
\textbf{Christopher~M Herald}, \emph{Flat connections, the {A}lexander
invariant, and {C}asson's invariant}, Comm. Anal. Geom. 5 (1997) 93--120

\bibitem{herald:perturbations}
\textbf{Christopher~M Herald}, \emph{{Legendrian cobordism and Chern--Simons
theory on $3$--manifolds with boundary}}, Comm. Anal. Geom. 2 (1994) 337--413

\bibitem{herald:pa9}
\textbf{Christopher~M Herald}, \emph{Transversality for equivariant gradient
systems and gauge theory on 3--manifolds} (2003), preprint

\bibitem{levine:polynomials}
\textbf{Jerome~P Levine}, \emph{Polynomial invariants of knots of codimension
two}, Ann. of Math. (2) 84 (1966) 537--554

\bibitem{masataka}
\textbf{Kaneko Masataka}, \emph{Casson's knot invariant and gauge theory},
Topology Appl. 112 (2001) 111--135

\bibitem{neumann-wahl}
\textbf{Walter~D Neumann}, \textbf{Jonathan Wahl}, \emph{Casson invariant of
links of singularities}, Comment. Math. Helv. 65 (1990) 58--78

\bibitem{neumann:plumbing}
\textbf{Walter~D Neumann}, \emph{An invariant of plumbed homology spheres},
from: ``Topology Symposium, Siegen 1979 (Proc. Sympos., Univ. Siegen, Siegen,
1979)'', Lecture Notes in Math. 788, Springer, Berlin (1980)  125--144

\bibitem{nicolaescu:swbook}
\textbf{Liviu~I Nicolaescu}, \emph{Notes on {S}eiberg-{W}itten theory},
   American Mathematical Society, Providence, RI (2000)

\bibitem{ruberman:ds}
\textbf{Daniel Ruberman}, \emph{Doubly slice knots and the {Casson--Gordon}
invariants}, Trans. Amer. Math. Soc. 279 (1983) 569--588

\bibitem{ruberman-saveliev:casson}
\textbf{Daniel Ruberman}, \textbf{Nikolai Saveliev}, \emph{Rohlin's invariant
and gauge theory I. Homology $3$-tori} (2003), {Comm. Math. Helv.} (to
appear). {\tt http://front.math.ucdavis.edu/math.GT/0302131}

\bibitem{salamon:lefschetz}
\textbf{Dietmar~A Salamon}, \emph{Seiberg-{W}itten invariants of mapping tori,
symplectic fixed points, and {L}efschetz numbers}, from: ``Proceedings of
6th G\"okova Geometry-Topology Conference'', Turkish J. Math. 23 (1999)
117--143

\bibitem{saveliev:brieskorn}
\textbf{Nikolai Saveliev}, \emph{Floer homology of {B}rieskorn homology
spheres}, J. Differential Geom. 53 (1999) 15--87

\bibitem{saveliev:seifert}
\textbf{Nikolai Saveliev}, \emph{Representation spaces of {S}eifert fibered
homology spheres}, Topology Appl. 126 (2002) 49--61

\bibitem{saveliev:cork}
\textbf{Nikolai Saveliev}, \emph{A note on {A}kbulut corks}, Math. Res. Lett.
10 (2003) 777--786

\bibitem{siebenmann:rohlin}
\textbf{Laurent Siebenmann}, \emph{On vanishing of the Rohlin invariant and
non-finitely amphicheiral homology 3-spheres}, from: ``Topology Symposium,
Siegen 1979 (Proc. Sympos., Univ. Siegen, Siegen, 1979)'', Lecture Notes in
Math. 788, Springer, Berlin (1980)  172--222

\bibitem{taubes:casson}
\textbf{Clifford~H Taubes}, \emph{Casson's invariant and gauge theory},
J. Differential Geom. 31 (1990) 547--599

\bibitem{wang:involution}
\textbf{Shuguang Wang}, \emph{Moduli spaces over manifolds with involutions},
Math. Ann. 296 (1993) 119--138

\end{thebibliography}

\end{document}